\documentclass[a4paper,12pt,twoside]{amsart}
\usepackage{amsmath,amsthm,amsfonts}
\usepackage[T1]{fontenc}
\usepackage{graphics,graphicx}
\usepackage{epsfig,color}

\newtheorem{thm}{Theorem}[section]

\newtheorem{cor}[thm]{Corollary}
\newtheorem{lem}[thm]{Lemma}
\newtheorem{lemma}[thm]{Lemma}
\newtheorem{prop}[thm]{Proposition}

\newtheorem{hypo0}{}

\theoremstyle{definition}

\newcommand{\Cal}{\mathcal}
\newcommand{\cal} {\mathcal}

\newcommand{\R}{{\mathbb{R}}}

\newcommand{\E}{{\mathbb{E}}}
\newcommand{\Q}{{\mathbb{Q}}}

\newcommand{\Z}{{\mathbb{Z}}}
\newcommand{\N}{{\mathbb{N}}}
\newcommand{\T}{{\mathbb{T}}}
\newcommand{\PP}{{\mathbb{P}}}

\newcommand{\tP}{\tilde{P}}
\newcommand{\tT}{\tilde{T}}

\def\eop{\qed}
\def\Proof {\vskip -2mm {\it Proof}.}
\def\proof {\vskip -2mm {\it Proof}.}

\def \notdiv {\not \hskip 2pt | \,}

\title [Diffusive behaviour of ergodic sums over rotations] {Diffusive behaviour of ergodic sums over rotations}

\date {\today}
\author{Jean-Pierre Conze}
\author{Stefano Isola}
\author{St\'ephane Le Borgne}
\address{IRMAR, UMR CNRS 6625, Universit\'e de Rennes I
\vskip 0mm Campus de Beaulieu, 35042 Rennes Cedex, France}
\email{conze@univ-rennes1.fr}
\email{stephane.leborgne@univ-rennes1.fr}
\address{ Scuola di Scienze e Tecnologie, Universit\'a di Camerino,
\vskip 0mm  via Madonna delle Carceri, 62032 Camerino, Italy}
\email{stefano.isola@unicam.it}

\dedicatory{Dedicated to the memory of Eugene Gutkin}

\subjclass[2010]{Primary: 11A55, 42A55, 60F05, 60F17; Secondary: 37D50.} 
\keywords{rotation, subsequences, variance, central limit theorem, lacunary series, almost sure invariance principle, periodic rectangular billiard}

\begin{document}

\begin{abstract} For a rotation by an irrational $\alpha$ on the circle and a BV function $\varphi$, we study the variance of the ergodic
sums $S_L \varphi(x) := \sum_{j=0}^{L -1} \, \varphi(x + j\alpha)$. When $\alpha$ is not of constant type, we construct sequences $(L_N)$
such that, at some scale, the ergodic sums $S_{L_N} \varphi$  satisfy an ASIP. Explicit non-degenerate examples are given,
with an application to the rectangular periodic billiard in the plane.
\end{abstract}

\maketitle

\tableofcontents

\section*{\bf Introduction}

Given a measure preserving map $T$ on a probability space $(X, \cal A, \mu)$, many results link the stochasticity
of $T$ to limit theorems in distribution for the ergodic sums of an observable $\varphi$ on $X$. A simple example is
$T: x \to 2 x \mod 1$ on $X = \R/\Z$ endowed with the Lebesgue measure: the normalized ergodic sums satisfy a Central Limit Theorem (CLT)
when $\varphi$ is H\"older or with bounded variation.

When $T$ is an irrational rotation $ x \to x +\alpha {\rm \ mod \ }1$ on $X$, the picture is quite different.
Depending on the Diophantine properties of $\alpha$, too much regularity for $\varphi$ can imply that $\varphi$ is a coboundary.
In that case, there is no way to normalize its ergodic sums in order to get a non-degenerate asymptotic distribution. Therefore, it is natural to
consider less regular but BV (bounded variation) functions, in particular step functions. Nevertheless, by the Denjoy-Koksma inequality, the ergodic sums
$S_L\varphi(x) = \sum_0^{L-1} \varphi(x +j \alpha)$ of a BV function $\varphi$ are uniformly bounded along the sequence $(q_n)$ of denominators of
$\alpha$. This leads to consider other subsequences $(L_n)$ with the hope that, along  $(L_n)$, there is a diffusive behaviour at some scale
for the ergodic sums.

Results on the CLT in the context of Fourier series, which are related to our framework, trace back to Salem and Zygmund \cite{SaZy48} in the 40's.
M. Denker and R. Burton in 1987, then M. Weber, M. Lacey and other authors gave results on a CLT for ergodic sums generated by rotations.
But their goal was the construction of some functions, necessarily irregular, whose ergodic sums satisfy a CLT after self-normalization.
The limit theorem along subsequences that we show below is for simple steps functions. In this direction, for
$\psi := 1_{[0, \frac12[} - 1_{[\frac12, 0[}$, F. Huveneers \cite{Hu09} proved that for every $\alpha \in \R \setminus \Q$ there is a sequence
$(L_n)_{n \in \N}$ such that $S_{L_n}\psi / \sqrt n$ is asymptotically normally distributed.

Here we consider the diffusive behaviour of the ergodic sums of step functions $\varphi$ over a rotation by $\alpha$. We study growth
of the mean variance and approximation of subsequences at a certain scale by a Brownian motion when $\alpha$ is not of constant type.
Another method, including the constant type case, will be presented in a forthcoming paper.

The content of the paper is the following. An analysis of the variance along subsequences is given
in Section \ref{varestmSect}. Then, in Section \ref{approxLac}, we introduce an approximation by lacunary series, in order to show
in Section \ref{cltRot} an ASIP (almost sure invariance principle) for subsequences of ergodic sums for BV observables when $\alpha$
has unbounded partial quotients. The method relies on the stochastic behaviour of sums of the form $\sum_1^N f_n(k_n x)$ where $(k_n)$
is a fast growing sequence of integers and $(f_n)$ a bounded set of functions in a class which contains the BV functions. It is based on a result
of Berkes and Philipp \cite{BePh79} in a slightly extended version (see appendix, Section \ref{ASIP0}).

Examples with a non-degenerate limit in distribution are presented in Section \ref{appli1}. The result has an application to a geometric model,
the billiard flow in the plane with periodic rectangular obstacles when the flow is restricted to special directions.

To conclude this introduction, let us observe that result presented here for an isometric map, the rotation by $\alpha$,
is related to the dimension 1. For instance for rotations on the 2-torus, the natural framework is to consider, instead of the single rotation,
the $\Z^2$-action generated by two independent rotations.

\vskip 2cm
\goodbreak 
{\bf Preliminaries}

In what follows, $\alpha$ will be an irrational number in $]0, 1[$. Its continued fraction expansion is denoted by
$\alpha = [0; a_1, a_2, ..., a_n, ...]$. We will need some reminders about continued fractions (see, for instance, \cite{Kh37}).

For $u\in\R$, set $\|u\| := \inf_{n \in \Z} |u - n|$.
Let $(p_n/q_n)_{n \ge 0}$ be the sequence of the {\it convergents} of $\alpha$. The integers $p_n$ (resp. $q_n$) are the {\it numerators}
(resp. {\it denominators}) of $\alpha$. They satisfy the following relations: $p_{-1}=1$, $p_0=0$, $q_{-1}=0$, $q_0=1$ and
\begin{equation} \label{converg_eq}
q_{n+1} = a_{n+1} q_n + q_{n-1}, \ p_{n+1} = a_{n+1} p_n + p_{n-1}, \ (-1)^n = p_{n-1} q_n - p_n q_{n-1}, \ n \geq 0.
\end{equation}
Let $n \geq 0$. We have, $\|q_n \alpha\| = (-1)^n (q_n \alpha - p_n)$,  $1 = q_n\|q_{n+1} \alpha \| + q_{n+1} \|q_n \alpha\|$.

Hence, setting $\displaystyle \alpha ={p_n \over q_n} + {\theta_n \over q_n}$, it holds
\begin{eqnarray}
&& {1\over 2 q_{n+1}} \leq {1\over q_{n+1}+q_n} \le \|q_n \alpha\|  = |\theta_n| \le {1\over q_{n+1}}
= {1 \over a_{n+1} q_n+q_{n-1}}. \label{f_3}
\end{eqnarray}
Moreover for $0 \leq k < q_n$, in every interval $[{k \over q_n}, {k+1 \over q_n}[$, there is a unique point of the form
$j \alpha {\rm \ mod \ } 1$, with $j \in \{0, ..., q_n-1\}$ and we have
\begin{eqnarray}
{1\over 2 q_{n}} \leq \|q_{n-1} \alpha \| \leq \|k \alpha \|, \text{ for }\ 1\le k < q_{n}. \label{f_4}
\end{eqnarray}

Recall that $\alpha$ is of constant type (or has bounded partial quotient (bpq)), if $\sup_k a_k < \infty$.

The uniform measure on $\T^1$ identified with $X = [0, 1[$ is denoted by $\mu$. We will denote by $C$ a generic constant which may change
from a line to the other. The arguments of the functions are taken modulo 1.
For a 1-periodic function $\varphi$, we denote by $V(\varphi)$ the {\it variation} of $\varphi$ computed for its restriction to the interval $[0, 1[$
and use the shorthand BV for ``bounded variation''. An integrable function
$\varphi$ is {\it centered} if $\mu(\varphi)= 0$.

Let $\Cal C$ be the class of centered BV functions. If $\varphi$ is in $\Cal C$, its Fourier coefficients $c_{r}(\varphi)$ satisfy:
\begin{eqnarray}
c_{r}(\varphi) = {\gamma_r(\varphi) \over r}, {\rm \ with \ }  K(\varphi) :=  \sup_{r \not = 0} |\gamma_r(\varphi)|  < +\infty. \label{majC}
\end{eqnarray}

The class $\Cal C$ contains in particular the step functions with a finite number of discontinuities.
An example which satisfies (\ref{majC}), but is not BV, is the 1-periodic function $\varphi$ such that
$\varphi(x)= -\log |x|$, for $x \in [-\frac12, \frac12[$, $x \not = 0$.

\vskip 3mm
{\it Notation}

Let $\varphi$ be in $\cal C$. The ergodic sums $\sum_{j=0}^{N-1} \varphi (x+j\alpha)$ are denoted by $S_N\varphi(x)$ or $\varphi_N(x)$. Hence we have
\begin{eqnarray}
\varphi_\ell(x) &:=  \sum_{j=0}^{\ell-1} \varphi(x + j \alpha) =
\sum_{r \not = 0} {\gamma_r(\varphi)\over r} \, e^{\pi i (\ell-1) r \alpha} \, {\sin \pi \ell r \alpha \over \sin \pi r \alpha} \, e^{2\pi i r x}.
\label{ergSum}
\end{eqnarray}

Let $q$ be such that for $p/q$ a rational number in lowest terms, $\|\alpha - p/q\| < {1 / q^2}$ (in particular for $q$ a denominator of $\alpha$).
By {\it Denjoy-Koksma inequality} we have
\begin{eqnarray}
\|\varphi_q\|_\infty = \sup_x |\sum_{\ell = 0}^{q-1}\varphi(x+\ell \alpha)| \le V(\varphi). \label{f_8}
\end{eqnarray}
In the $L^2$ setting, for functions which satisfy (\ref{majC}),
we have $\|\varphi_{q_n}\|_2  \leq 2 \pi \, K(\varphi)$. (See the remark after Proposition \ref{approxqn}).

Therefore, the size of the ergodic sums $\varphi_\ell$ depends strongly on the values of $\ell$, since for a BV centered function, the ergodic sums
are uniformly bounded along the sequence $(q_n)$ of denominators of $\alpha$.

Let us recall the Ostrowski expansion which gives a bound for the growth of the ergodic sums of a BV function.

Let $N \geq 1$ and $m=m(N)$ be such that $N \in [q_m, \ q_{m+1}[$.

We can write $N = b_m q_m + r$, with $1 \leq b_m \leq a_{m+1}$, $0 \leq r < q_{m}$. By iteration, we get for $N$ the following representation:
\begin{eqnarray*}
N =\sum_{k=0}^m b_k \, q_{k}, \text{ with } 0 \leq b_0 \leq a_1 -1, \ 0 \leq b_k \leq a_{k+1} \text{ for } 1 \leq k < m,
\  1 \leq b_m \leq a_{m+1}.
\end{eqnarray*}
Therefore, the ergodic sum can be written:
\begin{eqnarray}
S_N\varphi(x)= \sum_{\ell=0}^m \, \sum_{j=N_{\ell-1}}^{N_{\ell} -1}\varphi (x+j\alpha)
= \sum_{\ell=0}^m \, \sum_{j=0}^{b_\ell \, q_{\ell}-1} \varphi (x+N_{\ell-1} \alpha + j\alpha), \label{decompSul}
\end{eqnarray}
with $N_0 = b_0$, $N_\ell =\sum_{k=0}^\ell b_k \, q_{k}$ for $\ell \leq m$. It follows, for every $x$:
\begin{eqnarray}
|\sum_{j=0}^{N-1} \varphi(x+j\alpha)|  \leq \sum_{k=0}^{m(N)} \, b_k \, \|\varphi_{q_k}\|_\infty \leq  V(\varphi) \, \sum_{k=1}^{m(N)}\, b_k. \label{OstroBnd1}
\end{eqnarray}

{\it The aim}

In view of the Ostrowski expansion, we can ask if there is a diffusive behaviour of the ergodic sums along suitable subsequences defined in terms
of the $q_n$'s. We will see that, for rotations of non constant type, there are such sequences along which a CLT holds and that even a stronger stochastic
behaviour can occur: after redefining $(S_{L_n})$ on a probability space, denoting by $(\zeta(t)_{t \geq 0})$
the standard 1-dimensional Wiener process, we will show the existence of a sequence $(L_n)$ such that, for a sequence of r.v. $(\tau_n)$
and a constant $\lambda > 0$, we have $S_{L_n} = \zeta(\tau_n) + o(n^{\frac12 - \lambda}), \text{ a.e.}$
with $\tau_n / \|S_{L_n}\|_2^2 \to 1$, a.e. ($\|S_{L_n}\|_2^2$ is of order $n$).

The method of proof is the following. As shown below, if $a_{n}$ is big, there is $f_n$ such that the
${1\over q_n}$-periodic function $f_n(q_n .)$ approximates well the ergodic sum $\varphi_{q_n}$, cf. (\ref{majapproxqn20}).
When $a_{n}$, or a subsequence of  $a_{n}$, is growing fast, a CLT for subsequences $(\varphi_{L_n})$ can be deduced from the stochastic properties
of sums of the form $\sum_{k=1}^ N f_{n_k} (q_{n_k}.)$. In the appendix, we will recall a result of Berkes and Philipp which provides a CLT
and an approximation by a Wiener process for sums of this form.

\section{\bf Variance estimates  for subsequences of ergodic sums \label{varestmSect}}

\vskip 3mm
\subsection{\bf Bounds for the variance} \label{boundVar}

\

We will use inequalities related to the repartition of the orbit of 0 under the rotation by $\alpha$.

\begin{lem} \label{majSums}
If $q_n$ is a denominator of $\alpha$ and $m \geq 1$, we have with an absolute constant $C$:
\begin{eqnarray}
\sum_{k: \, \|k \alpha \| \leq 1/m, \, k \geq q_n} {1 \over k^2} &\leq& C \, ({1 \over m q_n} + {1\over q_n^2}), \ \
\sum_{k: \,\|k \alpha \| \geq 1/m, \, k \geq q_n } {1 \over k^2} {1 \over \|k \alpha \|^2} \leq C \, ({ m \over q_n} + {m^2 \over q_n^2}),
\label{maj2} \\
\sum_{k=1}^{q_{n} -1} \ {1\over k^2} \, {1\over \| k \alpha\|^2} &\le& 6 \sum_{j=0}^{n-1} ({q_{j+1} \over q_j})^2.
\label{vari1}
\end{eqnarray}
\end{lem}
\proof \ Observe first that, if $f$ is  a nonnegative BV function with integral $\mu(f)$ and if $q$ is a denominator of $\alpha$, then:
\begin{eqnarray}
\sum_{k=q}^\infty {f(k\alpha) \over k^2} \leq {2\mu(f) \over q} + {2V(f) \over q^2}. \label{maj1}
\end{eqnarray}
Indeed, by (\ref{f_8}) applied to $f -\mu(f)$, $\sum_{k=q}^\infty {f(k\alpha) \over k^2}$ is less than
$$\sum_{j=1}^\infty  {1\over (j q)^2} \sum_{r=0}^{q-1} f((j q +r) \alpha)
\leq {1\over q^2} (\sum_{j=1}^\infty {1\over j^2})\, (q \, \mu(f) + V(f)) = {\pi^2 \over 6} ({\mu(f) \over q} + \, {V(f) \over q^2}).$$
Now, (\ref{maj2}) follows from (\ref{maj1}), taking $f(x)$ respectively $= 1_{[{0, {1\over m}]}}(|x|)$ and $= {1 \over x^2} 1_{ [{1\over m}, {1\over 2}[}(|x|)$.

For (\ref{vari1}) we can write the LHS as
\begin{eqnarray*}
&&\sum_{j=0}^{n-1} \, \sum_{\ell=0}^{q_{j+1} - q_j -1}  \
{1\over (q_j+\ell)^2} \, {1\over \| (q_j+\ell) \alpha\|^2}
\ \leq \ \sum_{j=0}^{n-1} \ {1\over q_j^2} \, \sum_{\ell=0}^{q_{j+1} -
q_j -1} \, {1\over \| (q_j+\ell) \alpha\|^2}.
\end{eqnarray*}
Using the fact that there is only one value of $r \alpha \hbox{ mod
}1$, for $1\le r < q_{j+1}$, in each interval $[{k \over q_{j+1}},
{k+1 \over q_{j+1}}[$, $k=1,...,q_{j+1}-1$, we have the following
bound which implies (\ref{vari1}):
\begin{eqnarray*}
\sum_{\ell=0}^{q_{j+1} - q_j -1} \, {1\over \| (q_j+\ell) \alpha\|^2} \le {1 \over \|q_j \alpha\|^2} + \sum_{k=1}^\infty {1
\over (k/q_{j+1})^2}  \leq (q_j+q_{j+1})^2 + {\pi^2 \over 6} q_{j+1}^2 \leq 6 q_{j+1}^2 . \eop
\end{eqnarray*}

Now, we study the behaviour of the variance for the ergodic sums $\varphi_n (x)=\sum_{j=0}^{n-1}\varphi(x+j\alpha)$
of a function $\varphi(x) = \sum_r \frac{\gamma_r }{r} e^{2\pi i r x}$ in the class $\cal C$. We have:
\begin{eqnarray}
\Vert \varphi_n\Vert_2^2 = \sum_{r }\frac{|\gamma_r |^2}{r^2} G_n(r\alpha), \quad \hbox{with} \quad G_n(t)
:= {\sin^2 n\, \pi \,t \over \sin^2   \pi \, t}. \label{Gnx0}
\end{eqnarray}
The $1$-periodic function $G_n$ satisfies $\int_0^1 G_n(t) \, dt = n$ and the symmetry
$G_n(t)=G_n(1-t)$ for $0\leq t \leq 1$, so that $G_n(\{r \alpha\})= G_n( \Vert r \alpha \Vert)$.
We set
$$\langle G_n(t)\rangle :={1\over n} \sum_{k=0}^{n-1} G_k (t)= {1\over \sin^2 \pi t}\left[ {1\over 2}-{1\over 4n}
\left( 1+{\sin (2n-1)\pi t \over \sin\pi t}\right)\right].$$
The mean satisfies the following lower bounds:
\begin{eqnarray}
&&\langle G_n(t)\rangle \geq \frac {n^2}{\pi^2}, \text{ for } 0\leq t \leq {1\over 2 n},
\ \ \geq {1\over 8\, \pi^2 t^2}, \text{ for } {1\over 2 n} \leq t \leq {1\over 2 }. \label{lb1}
\end{eqnarray}

\vskip 3mm \begin{lem} (upper bound)  There is a constant $C$ such that, if $\varphi$ satisfies (\ref{majC}),
\begin{equation} \|\varphi_n\|_2^2 \le C K(\varphi)^2 \, \sum_{j=0}^\ell a_{j+1}^2, \ \forall n \in [q_\ell ,  q_{\ell+1}[.
\label{varipsi}\end{equation}
\end{lem}
\Proof \ Using (\ref{maj2}), we have (with the last inequality satisfied if $q_\ell \leq n)$:
\begin{eqnarray*}
\frac12\sum_{|k| \geq q_\ell} {1 \over k^2} {\|n k\alpha\|^2 \over
\|k\alpha\|^2} &\leq&  n^2 \sum_{\|k \alpha\| \leq 1/n, \, k \geq
q_\ell} {1 \over k^2} +  \sum_{\|k \alpha\| > 1/n, \, k \geq q_\ell} {1 \over k^2} {1 \over \|k\alpha\|^2} \\
&\leq&  {n^2 \over n q_\ell} +  {n^2\over q_{\ell}^2} + { n \over
q_\ell} +  {n^2 \over q_\ell^2} =  2 {n \over q_\ell} + 2 {n^2 \over q_\ell^2} \leq 4 {n^2 \over q_\ell^2}.
\end{eqnarray*}
Let $\varphi$ satisfy (\ref{majC}). Let $q_\ell \le n < q_{\ell+1}$.
From (\ref{vari1}) and (\ref{maj2}) of Lemma \ref{majSums}, we have
\begin{eqnarray*}\|\varphi_n\|_2^2 &=& \sum_{k\not = 0} |c_{k}(\varphi)|^2
{|1 - e^{2\pi i n k \alpha}|^2 \over |1 - e^{2\pi i k \alpha}|^2}
\leq K(\varphi)^2 \sum_{k\not = 0} {1 \over k^2} {\|n k \alpha\|^2 \over \|k\alpha\|^2} \\
&\leq& K(\varphi)^2 \sum_{0 < k < q_{\ell}} {1 \over k^2} {1 \over
\|k\alpha\|^2} + K(\varphi)^2 \sum_{k \geq q_{\ell}} {1 \over k^2} {\|nk\alpha\|^2 \over \|k\alpha\|^2} \\
&\le& K(\varphi)^2 [\sum_{j=0}^{\ell-1} ({q_{j+1} \over q_j})^2] + 6 K(\varphi)^2 {n^2 \over q_{\ell}^2}
\leq C K(\varphi)^2 \sum_{j=0}^\ell a_{j+1}^2. \eop
\end{eqnarray*}
When $\alpha$ is of bounded type, this gives $\max_{q_\ell \le n < q_{\ell+1}} \|\varphi_n \|_2 = O(\ell)$.

\vskip 3mm
For the mean $\langle D\varphi\rangle_n$ of the square norm of the ergodic sums, we have by (\ref{Gnx0}):
\begin{equation} \label{media0}
\langle D\varphi\rangle_n :={1\over n} \sum_{k=0}^{n-1} \Vert \varphi_k\Vert_2^2
=\sum_{r \not = 0}\frac{|\gamma_r |^2}{r^2} \langle G_n( \Vert r \alpha \Vert)\rangle.
\end{equation}

\begin{thm} \label{lowbd}  \cite{Is06} (lower bound) There is a constant $c> 0$ such that
\begin{eqnarray}
\langle D\varphi\rangle_n \geq c \sum_{j=0}^{\ell-1}|\gamma_{q_{j}}|^2 a_{j+1}^2, \ \forall n \in [q_\ell, q_{\ell+1}[. \label{minVariance0}
\end{eqnarray}
\end{thm}
\proof \ From (\ref{media0}) and (\ref{lb1}), we get the estimate
\begin{equation} \label{stima}
\langle D\varphi\rangle_n \geq  \sum_{r\, : \;  \Vert\, r\alpha\,\Vert \geq \frac  1{2n}}
{|\gamma_r|^2\over  8 \pi^2  (r \cdot \Vert r \alpha \Vert)^{2}}+ \sum_{r\ne 0\, : \;  \Vert\, r\alpha\,\Vert < \frac  1{2n}}
\frac {n^2  |\gamma_r|^2 }{\pi^2 r^2}.
\end{equation}
If $n \geq q_\ell$, then (\ref{f_3}) implies $\|q_j \alpha\| > {1 \over 2 n}$, for $j=1, ..., \ell-1$. Therefore, using
$q_{j} \, \Vert q_j \alpha \Vert <a_{j+1}^{-1}$, we have:
\begin{eqnarray*}
\langle D\varphi\rangle_n \geq C_1\sum_{\Vert\, r\alpha\,\Vert \geq \frac  1{2n}}
{ |\gamma_r|^2\over  (r \cdot \Vert r \alpha \Vert)^{2}}\geq c_1\,\sum_{j=0}^{\ell-1}
{ |\gamma_{q_j}|^2\over (q_{j}\cdot  \Vert q_j \alpha \Vert )^{2}} \geq C_2 \,\sum_{j=0}^{\ell-1}  |\gamma_{q_{j}}|^2 a_{j+1}^2. \eop
\end{eqnarray*}
By (\ref{minVariance0}), if $v_\ell$ is an integer $< q_{\ell+1}$ such that $\|\varphi_{v_\ell}\|_2 = \max_{k < q_{\ell+1}} \, \|\varphi_k\|_2$,
then
$$\|\varphi_{v_\ell}\|_2^2 \geq c \,\sum_{k=0}^{\ell-1}  |\gamma_{q_{k}}|^2 a_{k+1}^2.$$
If the sequence $(\gamma_{q_k}(\varphi))$ is bounded from below, the variance of most of the ergodic sums is of order of the scale
given by the $a_k$'s. If $\alpha$ has bounded partial quotients, the average of the variance grows at a logarithmic rate.

Examples will be given in Section \ref{appli1}.
To complete the picture we discuss a lower bound valid for functions $\varphi$ whose Fourier coefficients satisfy a definite (lower)
bound.

We define $u_n$ by
\begin{equation}\label{defilab}
u_n := \min \{\, u \, : \, \Vert\, q_{u}\alpha\,\Vert\,< \, {1\over 2n}\, \}=\max \{\, u \, : \,
\Vert\, q_{u -1}\alpha\,\Vert\,\geq\, {1\over 2n}\, \},
\end{equation}
The sequence $(u_n)$ cannot grow faster than $\log n$. Since
\begin{equation}\label{condiz}
\Vert\, r\alpha\,\Vert\geq \Vert\, q_{u_n -1}\alpha\,\Vert \geq {1\over2 n}, \  \forall r < q_{u_n},
\end{equation}
the integer  $q_{u_n}$ can be interpreted as the first time the orbit $\{x + \ell \alpha\}_{\ell \geq 1}$ returns into a neighborhood of
size $1/n$ of $x$.

We say that $\alpha$ is of {\sl type} $\gamma$ if $1\leq \gamma = \sup\{s : \liminf_{r\to \infty} r^s \cdot \| r\, {\alpha} \| =0\}$.

If $\alpha$ is of type $\gamma$, then $\liminf_{n\to \infty} {\log q_{u_n}\over \log n} = {1\over \gamma}$ (cf. \cite{Is06}).

\begin{lemma} \label{beta} Let $\alpha$ be of type $\gamma$ and $\varphi\in L^2$ be such that $|\gamma_{r}|> c\, r^{1-\beta}$ for some $\beta \in ]\frac 12, \gamma[$.
Then there is an infinite subset $I\subseteq \N$ and a constant $C>0$ such that
\begin{eqnarray}
\langle D\varphi\rangle_n  \geq C\, n^{ 2(1-\frac \beta {\gamma- \epsilon})(1+\beta -\frac \beta {\gamma- \epsilon})^{-1}}, \  n\in I, \  \forall \epsilon >0. \label{minDphi}
\end{eqnarray}
\end{lemma}
{\sl Proof.} We start again from the estimate (\ref{stima}).
Since $2q_{k} > 1/\Vert q_{k-1} \alpha \Vert > q_{k}$ for all $k$, using (\ref{defilab}) and (\ref{condiz}) we write
$$\langle D\varphi\rangle_n \geq  {1 \over  8 \pi^2 } \biggl( |\gamma_{q_{u_n-1}}|^2\, q_{u_n-1}^{-2}\, q_{u_n}^2 +
n^2 |\gamma_{q_{u_n}}|^2\, q_{u_n}^{-2}  \biggr).$$
The assumption on the Fourier coefficients of $\varphi$ then yields
$$\langle D\varphi\rangle_n  \geq  {c^2 \over  8 \pi^2 } \bigl( q_{u_n-1}^{-2\beta} \, q_{u_n}^2 + n^2  q_{u_n}^{-2\beta} \bigr).$$
On the other hand, one checks that the function
$$F_{a,b,s}(x):= x^{-2s} \bigl(a+ b x^2  \bigr), \ a,b>0,\ s \in (0,1), \  x\in \R^+,$$
satisfies
$$F_{a,b,s}(x) \geq \frac 1 {s^{s }(1-s)^{1- s}} \, a^{1-s} \, b^{s}, \ x\in\R^+,$$
and therefore
$$n^{-2s} \langle D\varphi\rangle_n  \geq \frac {c^2} { {8\pi^2} s^{s}(1-s)^{1-s}} q_{u_n-1}^{-2\beta (1-s)} q_{u_n}^{2(1-s-s\beta)}.$$
Now, if $\alpha$ is of type $\gamma \geq 1$, we have $\liminf_{r\to \infty} r^{\gamma -\epsilon} \cdot \Vert r \alpha\Vert=0$ for all $\epsilon >0$.
This implies that  $q_{{u_n}_j-1} \leq c_1 \, q_{{u_n}_j}^{\frac 1 {\gamma - \epsilon}}$ along an infinite subsequence $\{n_j\}$ and for some positive
constant $c_1$. Thereby, for each $s\in (0,1)$, we can find a constant $C_s$ so that the above yields:
\begin{eqnarray*}
&n_j^{-2s} \langle D\varphi\rangle_{n_j} \geq C_s\, q_{u_{n_j}}^{2\left(1-s-s\beta-\frac {\beta (1-s)}{\gamma - \epsilon} \right)}.
\end{eqnarray*}
Hence, taking $s$ such that the exponent at right is zero, we get (\ref{minDphi}). \eop

{\it Remark}:  Since $\varphi \notin {\cal C}$ for $\beta <1$, in order to apply Lemma \ref{beta} to some $\varphi \in \cal C$,
$\alpha$ must be of type $\gamma >1$.

\vskip 3mm
\subsection{\bf Lacunary series and approximation by ${1\over q_n}$-periodic functions} \label{approxLac}

\

We consider now another method giving information on the diffusive behaviour of the ergodic sums $S_N \varphi$.
It is based on approximation of $\varphi_{q_n}$ by ${1\over q_n}-$periodic functions.

Given $\varphi(x) = \sum_{r \not = 0} {\gamma_{r} \over r}\, e^{2\pi i r x}$, recall that the ergodic sums can be written:
\begin{eqnarray*}
\varphi_\ell(x) &= \sum_{r \not = 0} {\gamma_r(\varphi)\over r} \,  e^{\pi i (\ell-1) r \alpha}
\, {\sin \pi \ell r \alpha \over \sin \pi r \alpha} \, e^{2\pi i r x}.
\end{eqnarray*}
We use the following notations for $\ell \geq 1$:
\begin{eqnarray}
&&\widetilde \varphi_\ell(x) = \sum_{j=0}^{\ell-1} \varphi(x+\frac j\ell) = \sum_{r \not = 0} {\gamma_{r\ell} \over r}\, e^{2\pi i r\ell x}
= \widehat \varphi_{\ell}(\ell x), \text{ with }
\widehat \varphi_{\ell}(x) := \sum_{r \not = 0} {\gamma_{r\ell} \over r}\, e^{2\pi i rx}. \label{phihat}
\end{eqnarray}
In particular, $\widehat \varphi_{q_n}(x) := \sum_{r \not = 0} {\gamma_{r q_n} \over r}\, e^{2\pi i rx}$.
We will show that $\widehat \varphi_{q_n}(q_n.)$  is a ${1\over q_n}-$periodic approximation of $\varphi_{q_n}$, if $a_{n+1}$ is big.

{\it Remark}: If $\varphi \in \Cal C$, then $\widehat\varphi_{\ell}$ is also in $\cal C$ and satisfies: $K(\widehat\varphi_{\ell}) \leq K(\varphi)$.
If $\varphi$ is a BV function, then for every $\ell \geq 1$, the periodic function
$\widetilde \varphi_\ell = \widehat \varphi_{\ell}(\ell .)$ has the same variation on an interval of period and the variation of
$\widehat \varphi_\ell$ on $[0, 1[$ is less than $V(\varphi)$. When $\varphi$ has zero integral, this implies
\begin{eqnarray}
\|\widehat \varphi_\ell\|_\infty \leq V(\varphi). \label{boundper1}
\end{eqnarray}

\begin{prop} \label{approxqn} If $\varphi$  satisfies $(\ref{majC})$, then we have
\begin{eqnarray}
\|\varphi_{q_n} - \widehat \varphi_{q_n}(q_n.)\|_2^2 = \sum_{r \not = 0, \, q_n \notdiv r}  {1\over r^2}
\ {\|q_n r \alpha\|^2 \over \|r \alpha\|^2} = O(a_{n+1}^{-1}). \label{majapproxqn20}
\end{eqnarray}
\end{prop}
\proof \  If $\varphi$  satisfies $(\ref{majC})$, we have:
\begin{equation*}
\varphi_{q_n}(x) - \widehat \varphi_{q_n}(q_n x) = \sum_{r \not = 0} {\gamma_{r}(\varphi) \over r} \,  e^{\pi i (q_n-1) r \alpha}
\, {\sin \pi q_n r \alpha \over \sin \pi r \alpha} \, e^{2\pi i r x} - \sum_{r \not = 0} {\gamma_{q_n r}(\varphi) \over r}\, e^{2\pi i q_n rx}
= (A) + (B),
\end{equation*}
with
\begin{eqnarray*}
&&(A) = \sum_{r \not = 0} {\gamma_{q_n r}(\varphi) \over r} \,  [e^{\pi i (q_n-1) q_n r \alpha}
\, {\sin \pi q_n^2 r \alpha \over q_n \sin \pi q_n r \alpha}  - 1] \, e^{2\pi i q_n r x},  \\
&&(B) = \sum_{r \not = 0, q_n \not | r} {\gamma_{r}(\varphi) \over r} \,  e^{\pi i (q_n-1) r \alpha}
\, {\sin \pi q_n r \alpha \over \sin \pi r \alpha} \, e^{2\pi i r x}.
\end{eqnarray*}

Therefore, $\|\varphi_{q_n} - \widehat \varphi_{q_n}(q_n .)\|_2^2$ is equal to
\begin{eqnarray}
&&\sum_{r \not = 0} {|\gamma_{q_n r}(\varphi)|^2 \over r^2} \, |e^{\pi i (q_n-1) q_n r \alpha}
\, {\sin \pi q_n^2 r \alpha \over q_n \sin \pi q_n r \alpha}  - 1|^2 \,
+ \sum_{r \not = 0, q_n \not | r} {|\gamma_{r}(\varphi)|^2 \over r^2} \, |{\sin \pi q_n r \alpha \over \sin \pi r \alpha}|^2 \label{quadErrB}\\
&&\leq K(\varphi)^2\ \sum_{r \not = 0} {1 \over r^2} \, |e^{\pi i (q_n-1) q_n r \alpha}
\, {\sin \pi q_n^2 r \alpha \over q_n \sin \pi q_n r \alpha}  - 1|^2 \,
+ K(\varphi)^2 \, \sum_{r \not = 0, \, q_n \not | r} {1 \over r^2} \, |{\sin \pi q_n r \alpha \over \sin \pi r \alpha}|^2. \label{quadErrC}
\end{eqnarray}

Let $\varphi^0(x) = \{x\} -\frac12$. The Fourier series of $\varphi^0$ is
$\varphi^0(x) = {-1 \over 2\pi i} \, \sum_{r \not = 0} {1\over r} \  e^{2\pi i r x}$.

For $\delta_j \in [0, 1[$, we have: $|\varphi^0(x + \delta_j) - \varphi^0(x)| \leq \delta_j + 1_{[0, \delta_j]}(x)$; hence
$$\sum_{j=0}^{q_n-1} |\varphi^0(x + j \alpha +\delta_j) - \varphi^0(x + j \alpha)| \leq \sum_{j=0}^{q_n-1} \delta_j
+ \sum_{j=0}^{q_n-1} 1_{[0, \delta_j]}(x + j \alpha).$$
Since $|j \alpha - jp_n/q_n| \leq {1 \over a_{n+1} q_n}$, for $0 \leq j < q_n$, this implies:
$$|\varphi_{q_n}^0(x) - \widetilde \varphi_{q_n}^0(x)| \leq {1 \over a_{n+1}} + \sum_{j=0}^{q_n-1} 1_{[0, {1 \over a_{n+1} q_n}]}(x + j \alpha),$$
and therefore $\|\varphi_{q_n}^0 - \widehat \varphi_{q_n}^0(q_n.)\|_1 = 2 a_{n+1}^{-1}$.

As $\varphi_{q_n}^0$ and $\widetilde \varphi_{q_n}^0$ are bounded by $V(\varphi^0)$ (Denjoy-Koksma inequality and (\ref{boundper1})),
if follows:
\begin{eqnarray}
\|\varphi_{q_n}^0 - \widehat \varphi_{q_n}^0(q_n.)\|_2^2 \leq [\|\varphi_{q_n}^0\|_\infty + \|\widehat \varphi_{q_n}^0(q_n.)\|_\infty]
\ \|\varphi_{q_n}^0 - \widehat \varphi_{q_n}^0(q_n.)\|_1 \leq {4 \over a_{n+1}}. \label{majphi0}
\end{eqnarray}

By (\ref{quadErrB}) applied to $\varphi^0$, we obtain:
\begin{equation*}
\|\varphi^0_{q_n} - \widehat \varphi^0_{q_n}(q_n .)\|_2^2 = {1 \over 4 \pi^2} \, \sum_{r \not = 0} {1 \over r^2} \, |e^{\pi i (q_n-1) q_n r \alpha}
\, {\sin \pi q_n^2 r \alpha \over q_n \sin \pi q_n r \alpha}  - 1|^2 \,
+ {1 \over 4 \pi^2} \, \sum_{r \not = 0, \, q_n \not | r} {1 \over r^2} \, |{\sin \pi q_n r \alpha \over \sin \pi r \alpha}|^2.
\end{equation*}
It follows, by (\ref{quadErrC}) and (\ref{majphi0}):
\begin{eqnarray*}
&&\|\varphi_{q_n} - \widehat \varphi_{q_n}(q_n .)\|_2^2 \leq (2 \pi \, K(\varphi))^2 \, \|\varphi^0_{q_n} - \widehat \varphi^0_{q_n}(q_n .)\|_2^2
\leq ( 4 \pi \, K(\varphi))^2 \, a_{n+1}^{-1}. \eop
\end{eqnarray*}
{\it Remark}: Likewise, if $\varphi$  satisfies $(\ref{majC})$, then, since  $\|\varphi^0_{q_n}\|_2 \leq \|\varphi^0_{q_n}\|_\infty
\leq V(\varphi^0) =1$, we have
\begin{eqnarray*}
\|\varphi_{q_n}\|_2^2 &&= \sum_{r \not = 0} {|\gamma_{r}(\varphi)|^2 \over r^2} \, |{\sin \pi q_n r \alpha \over \sin \pi r \alpha}|^2
\leq K(\varphi)^2 \, \sum_{r \not = 0} {1 \over r^2} \, |{\sin \pi q_n r \alpha \over \sin \pi r \alpha}|^2\\
&&=  (2 \pi \, K(\varphi))^2 \, \|\varphi^0_{q_n}\|_2^2 \leq (2 \pi \, K(\varphi))^2 \,  V(\varphi^0)^2.
\end{eqnarray*}

\section{\bf CLT for rotations} \label{cltRot}

\subsection{\bf CLT and ASIP along subsequences for rotations}

\

Let  $(t_k)$ be a strictly increasing sequence of positive integers and let $(L_n)$ be the sequence of times defined by $L_0 = 0$ 
and $L_n = \sum_{k=1}^n q_{t_k}, \ n \geq 1$. Our goal is to show that, under a condition of the growth of $(a_n)$, the
distribution of $\varphi_{L_n}$ can be approximated by a Brownian motion. More precisely we will show the following ASIP (almost sure invariance principle, 
cf. \cite{PhSt75}) for $(\varphi_{L_n})$, when $a_{t_k +1}$ is fast enough growing:
\begin{thm} \label{TCLsubseq0} Let $(t_k)$ be a strictly increasing sequence of positive integers and  let $L_n = \sum_{k=1}^n q_{t_k}, n \geq 1$.
Assume the growth condition: $a_{t_k +1} \geq k^\beta$, with $\beta > 1$.
Then, for every $\varphi$ in the class $\cal C$ such that
\begin{eqnarray}
&&\sum_{k=1}^n \|\widehat \varphi_{q_{t_k}}\|_2^2 \geq c \, n, {\rm \ for \ a \ constant \ } c > 0, \label{min20}
\end{eqnarray}
we have $\|\varphi_{L_n} \|_2^2 / \sum_{k=1}^n \|\widehat \varphi_{q_{t_k}}\|_2^2 \to 1$ and the convergence in distribution (CLT)
\begin{eqnarray}
\varphi_{L_n}  /\|\varphi_{L_n} \|_2  \to \Cal N(0,1). \label{TCL0}
\end{eqnarray}

Moreover, keeping its distribution, the process $(\varphi_{L_n})_{n \geq 1}$
can be redefined on a probability space together with a Wiener process $\zeta(t)$ such that
\begin{eqnarray}
\varphi_{L_n} = \zeta(\tau_n) +o(n^{\frac12-\lambda}) \ \text{a.e.},
\end{eqnarray}
where $\lambda >0$ is an absolute constant and $\tau_n$ is an increasing sequence of random
variables such that $\tau_n/\|\varphi_{L_n} \|_2^2 \to 1$ a.e.
\end{thm}
\proof \ As in Ostrowski's expansion, the sum $\sum_{j=0}^{L_n-1} \varphi(x+j\alpha)$ reads
\begin{eqnarray*}
&& \sum_{k=0}^n \sum_{j=L_{k-1}}^{q_{t_k}+L_{k-1}-1} \varphi(x+j\alpha) =
\sum_{k=0}^n \sum_{j=0}^{q_{t_k}-1} \varphi(x+L_{k-1} \alpha + j \alpha) = \sum_{k=0}^n \varphi_{q_{t_k}}(x+L_{k-1} \alpha).
\end{eqnarray*}
We use Proposition \ref{approxqn} .
Let $g_k := |\sum_{j=L_{k-1}}^{q_{t_k}+L_{k-1}-1} \varphi(.+j\alpha) -  \widehat \varphi_{q_{t_k}}(q_{t_k}. + L_{k-1} \alpha)|$.
We have $\|g_k\|_2^2 \leq C a_{t_k+1}^{-1}$.

Assume that $a_{t_k +1} \geq k^\beta$. Then by Lemma \ref{boundgn} below, we have:
$\sum_{k=1}^n g_k =  O(n^{1 - {\beta \over 2}+ \varepsilon}) \text{ a.e.}$ Since $\beta > 1$,
this bound is comparable to the term of approximation in Theorem \ref{ext-BerkPhil} (appendix), which we apply now with 
$f_k = \widehat \varphi_{t_k}(. + L_{k-1} \alpha)$. Let us check the hypotheses of this theorem.

Condition (\ref{min20}) implies Condition (\ref{minMN}) of Theorem \ref{ext-BerkPhil} for $f_k$ and $n_k = q_{t_k}$
(see in the appendix Lemma \ref{varianceoN} which implies $\|\varphi_{L_n} \|_2^2 / \sum_{k=1}^n \|\widehat \varphi_{q_{t_k}}\|_2^2 \to 1$).

For (H\ref{hypofk}), observe that $\sup_k \|\widehat \varphi_{q_{t_k}} \|_\infty \leq V(\varphi) < +\infty$, $\sup_k \|\widehat \varphi_{q_{t_k}}\|_2
< +\infty$. Moreover, $|\gamma_{r q_{t_k}}| \leq K(\varphi)$  by (\ref{phihat}) and  there is a finite constant $C$ such that the tail of the Fourier series satisfies:
\begin{eqnarray*}
&&R(\widehat \varphi_{q_{t_k}}, t) \leq R(t), \ \forall k, \text { with } R(t) \leq C_R \, t^{-\frac12}.
\end{eqnarray*}
 For the sequence $(q_{t_k})$, the lacunarity condition (H\ref{hyponk}) as well as the arithmetic condition (H\ref{Gapo}) are satisfied,
 since $q_{t_{k+1}} / q_{t_k} > a_{t_k +1} \to \infty$.

Therefore, the hypotheses of Theorem \ref{ext-BerkPhil} are satisfied. The convergence in distribution (\ref{TCL0}) is a corollary. \eop

\begin{lem} \label{boundgn} Let $(g_n)$ be a sequence of nonnegative functions such that $\|g_n\|_2^2 = O(n^{-\delta})$, $\delta > 0$.
Then we have, for all $\varepsilon > 0$:
$$\sum_{k=1}^N g_k =  O(N^{1 - {\delta \over 2}+ {\varepsilon}}) \, a.e.$$
\end{lem}
\proof  \
For $\delta_1 = \frac12 - {\delta \over 2} + \varepsilon$, with $\varepsilon > 0$, we have convergence
$\int \, \sum_{k=1}^\infty \, {g_k^2 \over k^{2\delta_1}} \, d\mu < \infty$. Therefore,  $\sum_{k=1}^\infty {g_k^2 \over k^{2\delta_1}} =O(1)$, a.e.
which implies:
$$\sum_1^N g_k = \sum_1^N k^{\delta_1} {g_k \over k^{\delta_1}} \leq (\sum_1^N k^{2\delta_1})^\frac12 \,
(\sum_{k=1}^N {g_k^2 \over k^{2\delta_1}})^\frac12 =  O(N^{1 - {\delta \over 2}+ {\varepsilon}}), \, a.e. \eop$$

If the partial quotients of $\alpha$ are not bounded, then there is a sequence $(t_k)$ of positive integers tending to $+\infty$ such that
$a_{t_k +1} \geq k^\beta$, with $\beta > 1$. It follows:
\begin{cor} \label{coro1} Let $\alpha$ be an irrational rotation with unbounded partial quotients. Then there are $\lambda > 0$ and an increasing sequence of integers 
$(t_k)_{k \geq 1}$ such that, for the sequence, $L_n = \sum_{k=1}^N q_{t_k}$, $n \geq 1$, under the non-degeneracy condition (\ref{min20}), we have
\begin{eqnarray}
\varphi_{L_n} = \zeta(\tau_n) +o(n^{\frac12-\lambda}) \ \text{a.e.} \label{approxBrwn1}
\end{eqnarray}
\end{cor}

{\it Remarks}: 1) It still remains the question of Condition (\ref{min20}) that we will check for explicit step functions. We will have to estimate:
\begin{eqnarray}
&&\sum_{k=0}^{\ell-1}|\gamma_{q_k}(\varphi)|^2 a_{k+1}^2 \text{ (for the lower bound of the mean variance)}, \label{Mvar1}\\
&&\frac1N\sum_{k=1}^N \|\widehat \varphi_{q_{t_k}}\|_2^2
= \frac1N\sum_{k=1}^N \, [\sum_{r \not = 0} \, {1\over r^2} |\gamma_{r q_{t_k}}|^2] \ (\text{for the ASIP}) \label{Asip1}.
\end{eqnarray}
To get the ASIP along a subsequence $L_n = \sum_1^n q_{t_k}$, we need an increasing sequence $(t_k)$ of integers such that
$a_{t_k +1} \geq k^\theta$, with $\theta > 1$, and $\liminf_n \frac1n\sum_{k=1}^n \, [\sum_{r \not = 0} \, {1\over r^2} |\gamma_{r q_{t_k}}|^2] > 0$.
For this latter condition, it suffices that $\liminf_n \frac1n\sum_{k=1}^n \, |\gamma_{q_{t_k}}|^2 > 0$.

2) The result of Theorem \ref{TCLsubseq0} is valid more generally for sequences of the form $L_n = \sum_{k=0}^n c_k \, q_{t_k}$, 
where $(c_n)$ is a bounded sequence of non negative integers.

3) Let $\alpha$ be of Liouville type. Then, under a non degeneracy condition which is checked in the examples below, the variance along subsequences is "in average" 
of the order of $\sum_1^N a_k^2$ as shown by Theorem \ref{lowbd}, whereas the variance for the subsequences described in this section is much smaller and grows linearly.

4) If $\alpha$ is such that $a_n$ is of order $n^\beta$ with $\beta > 1$, we find a sequence $(L_n)$ for which Theorem \ref{TCLsubseq0} holds with 
of growth at most $\exp(c \, n \, \ln n)$ for some $c > 0$. 

5) For $\alpha$ with bounded partial quotients, a different approach is necessary for the CLT. It is based, as suggested in \cite{Hu09}, on a decorrelation property between
 the ergodic sums  at time $q_n$ for BV functions. The details will be given in a forthcoming paper.

\subsection{\bf Application to step functions} \label{appli1}

\

\vskip 3mm
{\it Example 1.} $\varphi(x) =\varphi^0(x) = \{x\} - \frac12 = {-1 \over 2\pi i} \, \sum_{r \not = 0} {1\over r} \  e^{2\pi i r x}$.

Here the above formulas (\ref{Mvar1}), (\ref{Asip1}) reduce to:
${1\over 4\pi^2} \, \sum_{k=0}^{\ell-1} a_{k+1}^2, \ \frac1N\sum_{k=1}^N \|\widehat \varphi_{q_{t_k}}\|_2^2 = {\pi^2 \over 6}$.
One easily deduces a non-degenerate CLT and ASIP along subsequences for non-bpq rotations.

Now we consider steps functions. The non-degeneracy of the variance is related to the Diophantine properties (with respect to $\alpha$) of its
discontinuities. If $\varphi$ is a step function: $\varphi = \sum_{j \in J} v_j \, (1_{I_j} - \mu(I_j)),$ with $I_j= [u_j, w_{j}[$, its Fourier
coefficients are
$$c_r = \sum_{j \in J} v_j {e^{-2\pi i r w_j} - e^{-2\pi i r u_{j}} \over 2 \pi i r}
= \sum_{j \in J} {v_j\over \pi r} \, e^{-\pi i r (u_j+w_j)} \, \sin(\pi r (u_j - w_j)), \ r \not = 0.$$

The growth of the mean variance is bounded from below by
\begin{eqnarray}
q_\ell \leq n < q_{\ell+1} \Rightarrow \langle D\varphi\rangle_n \geq C \sum_{k=1}^{\ell} \, |\sum_{j \in J} v_j (e^{2\pi i q_{k}
(w_j - u_{j})} -1)|^2 a_{k+1}^2. \label{minVariance2}
\end{eqnarray}
One can  try to set conditions on the coefficients $v_j$ and the endpoints of the partition such that $|\gamma_{q_k}| \asymp 1$
when $k$ belongs to some subsequence $J\subset \N$.

For example, if the $a_k$'s are bounded, then it can be shown, with an argument of equirepartition as below,
that for a.e. choice of the parameters the lower bound $\langle D \varphi \rangle_n \geq c \ln n$ holds.

Now we consider different particular cases for generic or special values of the parameter.

\vskip 3mm
{\it Example 2.} \ $\varphi = \varphi(\beta, .) = 1_{[0, \beta[} - \beta = \sum_{r \not = 0} {1 \over \pi r} e^{-\pi i r \beta}
\, \sin(\pi r \beta) \ e^{2 \pi i r .}$.

Therefore $\gamma_r(\varphi)= {1\over \pi} e^{- \pi i r \beta} \sin \pi r \beta$ and Theorem \ref{lowbd} yields
$$\langle D\varphi\rangle_n  \geq C_1 \,{\sum_{k=0}^{\ell-1}} a_{k+1}^2  \, \sin^2{( \pi   q_k \beta)} \geq 4 C_1 \,{\sum_{k=0}^{\ell-1}}
a_{k+1}^2  \, \Vert  q_k \beta\Vert ^2\, , \quad \forall n \in [q_\ell,q_{\ell +1}[ \, .$$

For the mean variance, putting, for $\delta > 0$, $J_\delta:= \{k:\Vert q_{k} \beta \Vert \geq \delta\}$, we have
\begin{equation} \label{sst}
\langle D\varphi\rangle_n  \geq C \delta^{2} {\sum_{k\in J_\delta \cap \,[1,\ell[}} a_{k}^2 \, , \quad \forall n \in [q_\ell,q_{\ell +1}[ \, .
\end{equation}

For the CLT, we have $\widehat \varphi_{q_n} = \sum_{r \not = 0} {1 \over \pi r} e^{-\pi i r q_n \beta}
\, \sin(\pi r q_n \beta) \ e^{2 \pi i r .}$,
\begin{eqnarray}
\|\widehat \varphi_{q_n}\|_2^2 &=& {1 \over \pi^2} \sum_{r \not = 0} {|\sin(\pi r q_n \beta)|^2 \over r^2} \geq {1 \over \pi^2}\,|\sin(\pi q_n \beta)|^2.
\end{eqnarray}
Since $(q_k)$ is a strictly increasing sequence of integers, for almost every $\beta$ in $\T$, the sequence $(q_{t_k} \beta)$ is uniformly
distributed modulo 1 in $\T^1$. For a.e. $\beta$ we have:
$$\lim_N {1 \over N} \sum_{k=1}^N \|\widehat \varphi_{t_k}\|_2^2 = \lim_N {1 \over N} \sum_{k=1}^N {1 \over \pi^2} \sum_{r \not = 0}
{|\sin(\pi r q_{{t_k}} \beta)|^2 \over r^2} = \frac16.$$
Hence, Theorem \ref{TCLsubseq0} implies:
{\it If $\alpha$ is not of bounded type, there exists a sequence $(q_{t_k})$ of denominators of $\alpha$ such that, for the
subsequence $L_N = q_{t_1}+...+q_{t_N}$, for a.e. $\beta$,
\begin{eqnarray*}
{\sqrt 6 \over \sqrt N} \sum_{j=1}^{L_N} \varphi(\beta, . +j\alpha)
\underset{distribution}\to \cal N(0,1).
\end{eqnarray*}}
{\it A special case ($\beta = \frac12$}):
$\varphi := 1_{[0, {1 \over 2}[} - 1_{[{1 \over 2},1[} = \sum_{r} {2 \over \pi i (2r+1)} \, e^{2\pi i (2r+1) .}$.

In this case, we have $\gamma_{q_k} = 0, \text{ if } q_k \text{ is even}, \ =  {2 \over \pi i}, \text{ if } q_k \text{ is odd}$.
If $a_{n+1} \to \infty$ along a sequence such that $q_n$ is odd, then there is a sequence $(L_n)$ for which
\begin{eqnarray*}
{1 \over \sqrt N} \sum_{j=1}^{L_N} \varphi(\frac12, . +j\alpha)
\underset{distribution}\to \cal N(0,1).
\end{eqnarray*}

\vskip 3mm
{\it Remark}. Degeneracy can occur even for a cocycle which generates an ergodic skew product on $\T^1 \times \R$ (and therefore is not a measurable
coboundary). Let us consider $1_{[0, \beta[} - \beta$ and the so-called Ostrowski expansion of $\beta$:
$\beta=\sum_{n\geq 0}b_{n}q_{n}\alpha \textrm{ mod }1, \text{ with } b_n \in \Z$. Then it can be shown that
$\sum_{n\geq 0}\frac{|b_{n}|}{a_{n+1}} <\infty$ implies $\lim_k \|q_k \beta\| = 0$. If $\alpha$ is not bpq, there is an uncountable set of $\beta$'s
 satisfying the previous condition, but ergodicity of the cocycle holds if $\beta$ is not in the countable set $\Z \alpha + \Z$.

\vskip 3mm
{\it Example 3.} Let $\varphi$ be the step function: $\varphi = \varphi(\beta, \gamma, .)
= 1_{[0, \ \beta ]} - 1_{[\gamma , \ \beta +\gamma ]}$. The Fourier coefficients are
$c_r(\varphi) = {2i \over \pi} {1 \over r} e^{-\pi ir(\beta +\gamma )} \, \sin(\pi r \beta ) \ \sin(\pi r \gamma )$. We have
\begin{eqnarray*}
\|\widehat\varphi_{q_k}\|_2^2 &=& {4 \over \pi^2} \, \sum_{r \not = 0} {1 \over r^2} \ |\sin(\pi r q_k \beta )|^2 \ |\sin(\pi r q_k \gamma )|^2.
\end{eqnarray*}
As above, since $(q_k)$ is a strictly increasing sequence of integers, for almost every $(\beta ,\gamma )$ in $\T^2$, the sequence $(q_{t_k}
\beta , q_{t_k} \gamma )$ is uniformly distributed in $\T^2$. We have for a.e. $(\beta ,\gamma )$:
\begin{eqnarray*}
\lim_n {1\over n} \sum_{k=1}^n \|\widehat \varphi_{q_{t_k}}\|_2^2 &=& {4 \over \pi^2} \sum_{r \not = 0}  \lim_n {1\over n} \sum_{k=1}^n
{|\sin(\pi r q_{t_k} \beta )|^2 \ |\sin(\pi r q_{t_k} \gamma )|^2 \over r^2}
\\ &&= {4 \over \pi^2} \sum_{r \not = 0} \int \int {|\sin(\pi r y)|^2 \ |\sin(\pi r z)|^2 \over r^2} \ dy \, dz = \frac13.
\end{eqnarray*}

This computation and Theorem \ref{TCLsubseq0} imply the following corollary for $\varphi(\beta ,\gamma, .)$:

{\it If $\alpha$ is not of bounded type, there exists a sequence $(q_{t_k})$ of denominators of $\alpha$ such that, for the
subsequence $L_N = q_{t_1}+...+q_{t_N}$, for a.e. $(\beta,\gamma )$
$${\sqrt 3\over \sqrt N} \sum_{j=1}^{L_N} \varphi(\beta ,\gamma, . +j\alpha)
\underset{distribution}\to \cal N(0,1).$$}

\vskip 3mm
{\it Example 4.} Let us take $\gamma = \frac12$, i.e., $\varphi = \varphi(\beta, \frac12, .) = 1_{[0, \ \beta]} - 1_{[\frac12, \  \beta+ \frac12]}$.
We have:
$$\varphi(x) = {2 \over \pi} \, \sum_{r} \, e^{-\pi i (2r+1) \beta} \, {\sin(\pi (2r+1) \beta) \over (2r+1)} \, e^{2\pi i (2r+1) x}.$$

Hence: $|\gamma_{q_k}| \sim \|q_k \beta\|$ if $q_k$ is odd, else $=0$.
This example is like Example 2, excepted the restriction to odd values of the frequencies.

The lower bound for the mean variance (Theorem \ref{lowbd}) gives in this case:
\begin{eqnarray}
q_\ell \leq n < q_{\ell+1} \Rightarrow \langle D\varphi\rangle_n \geq C \sum_{0 \leq k \leq \ell-1, \, q_k {\rm \ odd}} \|q_{k}
\beta\|^2 \, a_{k+1}^2
\end{eqnarray}
and we have, if $(q_{t_k})$ is a sequence of odd denominators, for a.e. $\beta$:
$$\lim_N {1 \over N} \sum_{k=1}^N \|\widehat \varphi _{t_k}\|_2^2 = \lim_N {1 \over N} \sum_{k=1}^N {4 \over \pi^2} \sum_{r}
{|\sin(\pi (2r+1) q_{{t_k}} \beta)|^2 \over (2r+1)^2} = {2\over \pi^2} \sum_{r \in \Z}  {1\over (2r+1)^2} = \frac12.$$

\vskip 3mm
{\it Example 5.} (Vectorial cocycle) Now, in order to apply it to the periodic billiard, we consider the vectorial function
$\psi=(\varphi^1, \varphi^2)$, where $\varphi^1$ and $\varphi^2$ are functions as in the example 4) with parameters
$\beta_1=\frac{\alpha}{2}$, $\beta_2=\frac{1}{2}-\frac{\alpha}{2}$:
\begin{eqnarray*}
\varphi^1&&={1}_{[0,\frac{\alpha}{2}]}-{1}_{[\frac{1}{2},\frac{1}{2}+\frac{\alpha}{2}]}=\frac{2}{\pi}\sum_{r\in\Z}e^{-\pi i(2r+1)\frac{\alpha}{2}}
\frac{\sin(\pi (2r+1)\alpha/2)}{2r+1}e^{2\pi i (2r+1) .},\\
\varphi^2&&={1}_{[0,\frac{1}{2}-\frac{\alpha}{2}]}-{1}_{[\frac{1}{2},1-\frac{\alpha}{2}]}=\frac{-2i}{\pi}\sum_{r\in\Z}e^{\pi i(2r+1)\frac{\alpha}{2}}
\frac{\cos(\pi (2r+1)\alpha/2)}{2r+1}e^{2\pi i(2r+1) .},
\end{eqnarray*}
If $q_k$ is even both $\widehat{\varphi^1_{q_k}}$ and $\widehat{\varphi^2_{q_k}}$ are null. If $q_k$ is odd, we have
$$\| \widehat{\varphi^1_{q_k}}\|_2^2=\frac{4}{\pi^2}\sum_{r\in\Z}\frac{|\sin(\pi (2r+1)\frac{\alpha}{2})|^2}{(2r+1)^2},\ \| \widehat{\varphi^2_{q_k}}\|_2^2=\frac{4}{\pi^2}\sum_{r\in\Z}\frac{|\cos(\pi (2r+1)\frac{\alpha}{2})|^2}{(2r+1)^2}.$$
Let $q_k$ be odd. We use (\ref{f_3}). We have
$$\| q_k\beta_1\|= \| q_k \frac{\alpha}{2}\|=\|\frac{p_k}{2}+\frac{\theta_k}{2} \|, \text{ hence }
\left| \| q_k\beta_1\|-\|\frac{p_k}{2} \|\right|\leq \left|\frac{\theta_k}{2}\right|\leq \frac{1}{2q_{n+1}}.$$
This implies:
\begin{eqnarray*}
\| \widehat{\varphi^1_{q_k}}\|_2^2&=&O(\frac{1}{q_{k+1}}),\ {\rm if }\ p_k\ {\rm is \ even },\\
\| \widehat{\varphi^1_{q_k}}\|_2^2&=&\frac{4}{\pi^2}\sum_{r\in\Z}\frac{1}{(2r+1)^2}+O(\frac{1}{q_{k+1}})=1+O(\frac{1}{q_{k+1}}),\ {\rm if }
\ p_k\ {\rm is \ odd }.
\end{eqnarray*}
Similarly, we have
$$\| q_k\beta_2\|=\|\frac{q_k}{2}-\frac{p_k}{2}-\frac{\theta_k}{2} \|,
\text{ hence } \left| \| q_k\beta_2\|-\left(\|\frac{1}{2}+\frac{p_k}{2} \|\right|\right)\leq \left|\frac{\theta_k}{2}\right|\leq \frac{1}{2q_{n+1}},$$
hence:
\begin{eqnarray*}
\|\widehat{\phi^2_{q_k}}\|_2^2&=&O(\frac{1}{q_{k+1}}), \ {\rm if }\ p_k\ {\rm is \ odd,}
\ \ \|\widehat{\phi^2_{q_k}}\|_2^2 = 1+O(\frac{1}{q_{k+1}}), \ {\rm if }\ p_k\ {\rm is \ even.}
\end{eqnarray*}

\begin{lem} For almost every $\alpha$, for $\beta > 1$, there exists a sequence $(t_k)$ such that
\begin{eqnarray}
 q_{t_k} \text{ is odd}, \ p_{t_{2k}} \text{ is even},  \ p_{t_{2k+1}}\text{ is odd} \text{ and } a_{t_k+1} \geq k^{\beta}
 \text{ for all } k. \label{oddEven1}
\end{eqnarray}
\end{lem}
\proof \  Let us examine, for three consecutive terms, the configurations for $(p_n, q_n)$ modulo 2, i.e., the parities of $p_n, q_{n}$,
$p_{n+1}, q_{n+1}$, $p_{n+2}, q_{n+2}$. Suppose that $a_{n+2}$ is even. Then, using (\ref{converg_eq}) we see  that the only possible
configurations are:
\begin{eqnarray*}
&&[(0, 1), (1, 0), (1, 1)], \ \  [(0, 1), (1, 1), (1, 0)], \ \  [(1, 0), (0, 1), (1, 1)], \\
&& [(1, 1), (0, 1), (1, 0)], \ \  [(1, 0), (1, 1), (0, 1)], \ \  [(1, 1), (1, 0), (0, 1)].
\end{eqnarray*}
Taking either $p_n, q_{n}$, $p_{n+1}, q_{n+1}$, or  $p_n, q_{n}$, $p_{n+2}, q_{n+2}$, or $p_{n+1}, q_{n+1}$, $p_{n+2}, q_{n+2}$, it follows that we
can find among three consecutive convergents $(p,q), (p', q')$ the desired parities, i.e., $p$ even, $q$ odd, $p'$ odd, $q'$ odd.

Now let be given $A_1, A_2, A_3$ three positive integers. Using the ergodicity of the Gauss map $x \to \{1 / x\}$ in the class of the Lebesgue measure,
one easily shows that for a.e. $x$, there are infinitely many values of $n$ such that $a_n= A_1$, $a_{n+1}= A_2$, $a_{n+2}= A_3$.
By choosing successively for $A_1, A_2, A_3$ arbitrary big integers and $A_3$ even, and using the above analysis of possible configurations,
we get a sequence $(t_k)$ such that the condition on the parities as in (\ref{oddEven1}) is satisfied and $ \lim_k a_{t_k+1}=+\infty$. 
Now, taking a subsequence, still denoted $(t_k)$, we can insure that $a_{t_k+1} \geq k^{\beta}$, with $\beta > 1$ for all $k \geq 1$.
\eop
\begin{thm} \label{vectCLT1} Let $\alpha$ be an irrational number satisfying the (generic) condition (\ref{oddEven1}) holds for a sequence  $(t_k)$
and  the convergents of $\alpha$. Then, for $L_N=\sum_{1}^{N}q_{t_k}$, the sequence $(N^{-1/2}\psi_{L_N})_{N\geq 1}$ satisfies a 2-dimensional CLT
with a non-degenerate diagonal covariance matrix $\begin{pmatrix}\frac12 & 0 \\ 0 & \frac12 \end{pmatrix}$.
\end{thm}
\proof \ \rm It suffices to prove the CLT for $\psi^{u,v}=u\varphi^1+v\varphi^2$ where $(u,v)$ is an arbitrary fixed vector. By hypothesis
the sequence $(q_{t_k})$ is superlacunary. Therefore (Lemma \ref{varianceoN} in the appendix), we have
$$\|\sum_{k=1}^N \widehat{\psi^{u,v}_{q_{t_k}}}\|_2^2=\sum_{k=1}^N \|\widehat{\psi^{u,v}_{q_{t_k}}}\|_2^2+o(N).$$
The above computations show that, under Condition (\ref{oddEven1}),
\begin{eqnarray*}
\|\widehat{\psi^{u,v}_{q_{t_k}}}\|_2^2&=&v^2+O(\frac{1}{q_{t_k+1}}), \ {\rm if }\ k\ {\rm is \ even },
\ \ \|\widehat{\psi^{u,v}_{q_{t_k}}}\|_2^2 = u^2+O(\frac{1}{q_{t_k+1}}), \ {\rm if }\ k\ {\rm is \ odd }.
\end{eqnarray*}
Thus the hypotheses of Theorem \ref{TCLsubseq0} are fulfilled and $\left(N^{-1/2}\left(u\varphi^1_{L_N}+v\varphi^2_{L_N}\right)\right)$ converges
in distribution to the centered Gaussian law with variance $\frac{u^2+v^2}{2}$.
\eop

\vskip 3mm
{\it Remarks:} 1) Using the continued fraction expansion of $\alpha$, one can find rotations and subsequences $(L_n)$
with arbitrary values of the covariance matrix in the above CLT.

2) We have only considered the CLT for the vectorial cocycle, but not the ASIP, since the result of Theorem \ref{ext-BerkPhil}
is one dimensional. Further work should be done to obtain a vectorial ASIP (approximation by a two dimensional Brownian motion).
Let us mention that a method to get it, could be to adapt the method developed by S. Gou\"ezel in \cite{Go10}.

\vskip 3mm
\subsection{\bf Application to the periodic billiard in the plane}

\

{\it Description of the model}

We start with a brief description of the billiard flow in the plane with $\mathbb{Z}^2$-periodically distributed obstacles.
The flow acts on the set of configurations, where a configuration is a position in the complementary of the obstacles in $\R^2$ together with a unitary speed vector.
The flow is defined according to the usual rule: the ball (geometrically reduced to a point)  moves with
constant speed in straight line between two obstacles and obeys the laws of reflection when it hits the edge of an obstacle.

If the obstacles are strictly convex with regular boundary and positive curvature, the flow has a stochastic behaviour
with a rate of diffusion of order $\sqrt{N}$ at time $N$. The situation is quite different when the boundary of the obstacles
is piecewise flat. We consider here the billiard with $\Z^2$-periodically distributed rectangular obstacles (sometimes called
rectangular Lorenz gas, cf. \cite{HW80}). We give below a short description of this model (cf. \cite{CoGu12}) and show that for special directions
the previous results on ergodic sums over rotations apply.

Let $0< a, b< 1$. For $(m,n)\in\Z^2$, let $R_{(m,n)}(a,b)\subset\R^2$ be the $a\times b$ rectangle centered at $(m,n)$ whose sides are parallel
to the coordinate axes and let $R(a,b)=ABCD$ be the rectangle of same size in the unit torus. The associated billiard flow $(T^t)_{t\geq 0}$ can be
viewed as the geodesic flow on the polygonal surface $\tP(a,b) = \R^2 \setminus \bigcup_{(m,n)\in\Z^2}R_{(m,n)}(a,b)$.
It acts on the set $U\tilde P(a,b)$ of unit tangent vectors and preserves the Liouville measure $\mu$ on it.

The quotient surface $P(a,b) = \tilde P(a,b)/\Z^2$ is the unit torus with a rectangular hole. It is rational and its dihedral group is $R_2$,
the group generated by two orthogonal reflections with the angle $\pi/2$ between their axes. The flow $(U\tilde P(a,b),T^t,\mu)$ decomposes 
as a one-parameter family  of {\em directional billiard flows} $(\tilde Z_{\eta},\tT^t_{\eta},\tilde \mu_{\eta})$, where $\eta\in[0,\pi/2]$.

Let $\eta\in(0,\pi/2)$ and let $X_{\eta}$ be the space consisting of unit vectors pointing outward, whose base points belong to $ABCD$
and whose directions belong to the set $\{\pm\eta,\pi\pm\eta\}$. A natural Poincar\'e's section $(\tilde X_{\eta},\tilde \tau_{\eta},\tilde \nu_{\eta})$
of the conservative part of the flow is given by its restriction to the configurations corresponding to the boundary of the obstacles, at times
of impact with the obstacles. We obtain so the {\it billiard map}.

{\it Rational directions and small obstacles} \label{small_obst_sub}

A direction $\eta\in[0,\pi/2]$ is rational if $\tan\eta\in\Q$. Rational directions $\eta(p,q) = \arctan (q/p)$, also simply denoted $(p,q)$,
correspond to pairs $(p,q)\in\N$ with relatively prime $p,q$.

In what follows we fix $(p,q)$ and assume that the following "{\it small obstacles condition}" is satisfied: $qa+pb \leq 1$.

The inequality above is strict if and only if the directional billiard flow
$(\tilde Z_{(p,q)},\tT^t_{(p,q)},\tilde \mu_{(p,q)})$ has a set of positive measure of orbits that do not encounter obstacles.

We will now investigate the Poincar\'e map $\tau_{(p,q)}:X_{(p,q)} \to X_{(p,q)}$ induced by the billiard map. We identify $X_{(p,q)}$ with 2 copies
of the rectangle $ABCD$: one copy carries the outward pointing vectors in the direction $\eta$ or $\pi+\eta$, the other one the outward pointing vectors
in the direction $\pi-\eta$ or $2\pi-\eta$.

One can reduce the model to the case of a direction of flow with angle $\eta = \pi/4$ and with the small obstacles condition $a+b \leq 1$ (see \cite{CoGu12}). 
Without loss of generality, we will consider this case. Let $\displaystyle \alpha:=\frac{a}{a+b}$.

The square of the Poincar\'e map for the direction $\pi/4$ can then be represented as two copies of the skew product
defined on $\T^1 \times \Z^2$ by $(x, z) \to (x + \alpha, z + \Psi(x))$, where the displacement function $\Psi$ is given by:
\begin{equation}   \label{cocy_Psi}
\Psi(x) = \left \{ \begin{array}{clcr}
(0, 1) \text{ for } x \in ]0,{1 \over 2}-{\alpha\over 2}[,\ \ (1, 0) \text{ for } x \in ]{1 \over 2}-{\alpha\over 2},{1 \over 2}[,\\
(0, -1)  \text{ for  } x \in ]{1 \over 2},1-{\alpha\over 2}[, \ \ (-1, 0) \text{ for } x \in ]1-{\alpha\over 2},1[.
\end{array} \right .
\end{equation}
This vectorial function reads $\Psi = (\psi_1, \psi_2)$, with $\psi_1= \varphi_1(. + \frac12 - \frac\alpha2)  = 1_{[\frac12 - \frac\alpha2, \frac12[}
- 1_{1 - \frac\alpha2,1[}$, $\psi_2= \varphi_2  = 1_{[0, \frac12 - \frac\alpha2[} - 1_{[\frac12, 1 - \frac\alpha2[}$, where $\varphi_1, \varphi_2$
are defined in Subsection \ref{appli1}.

The billiard map generates the cocycle $S(n, \Psi)(x) = \sum_{j=0}^{n-1} \, \Psi(x + j \alpha)$ (the ergodic sums of the displacement function
over the rotation by $\alpha$), which gives the label in the $\Z^2$-plane of the cell containing the ball after $2n$ reflections on the obstacles.

Let $\psi(x)$ be the length of the path for the ball starting  from $x \in ABCD$ (identified with the circle $\T$) up to its collision
with a second obstacle and let $X^t(x)$ be the position at time $t$ of a ball starting from $x$. We put
\begin{equation}  \label{centerpsi}
c = \int_\T \, \psi(x) \, dx,  \ \psi^0 = \psi - c.
\end{equation} 

The ergodic sum $T_n(x) = \sum_{j=0}^{n-1} \psi(x+j\alpha)$  is the hitting time of the ball with an obstacle after $2n$ collisions. 
Therefore $X^{T_n(x)} (x)$ belongs to the cell $S(n, \Psi)(x)$ in the $\Z^2$-plane.
 
\begin{figure}
\begin{center}
\includegraphics[scale=.3]{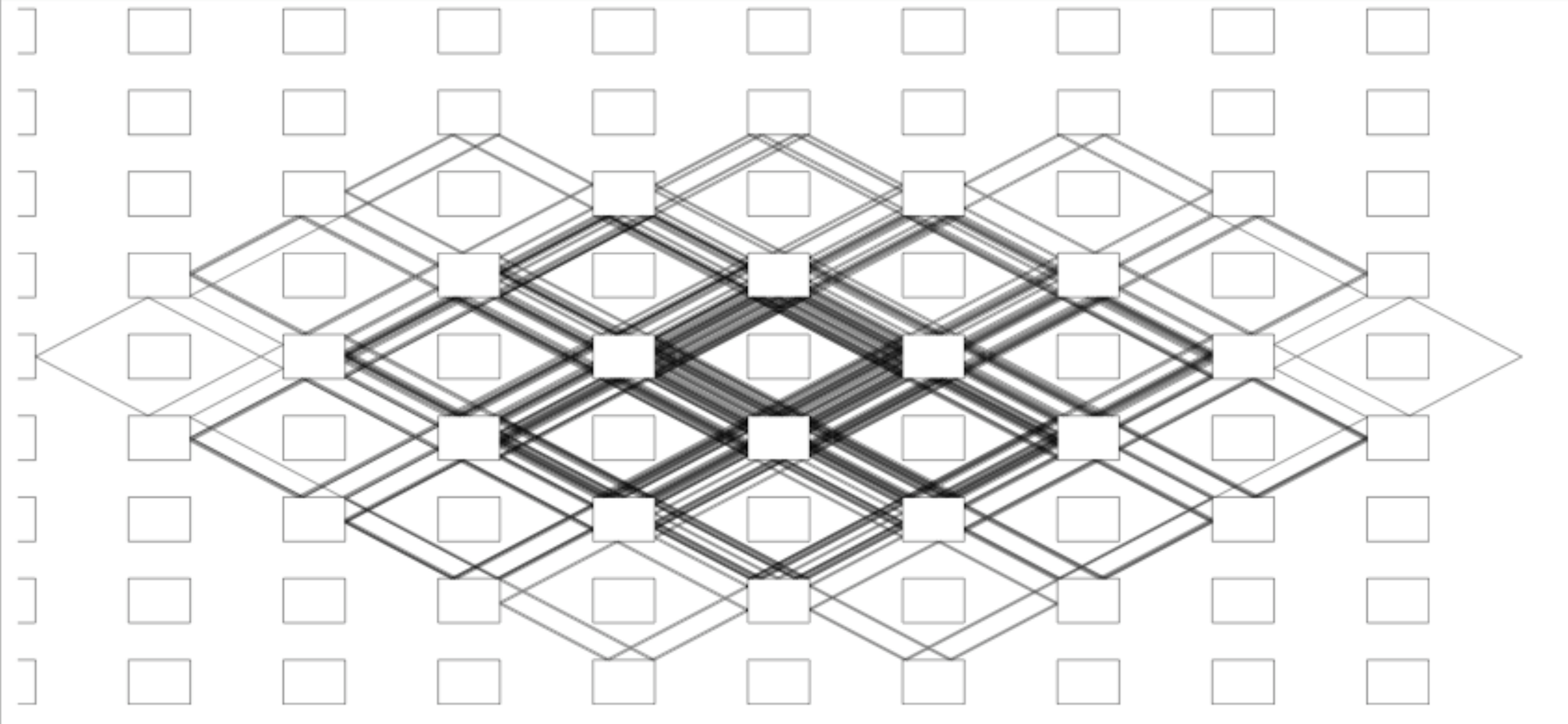}
\end{center}
\caption{\it An orbit of the rectangular billiard, angle $\pi /4$}
\label{bill-fig1}
\end{figure}

\begin{figure}
\begin{center}
\includegraphics[scale=.3]{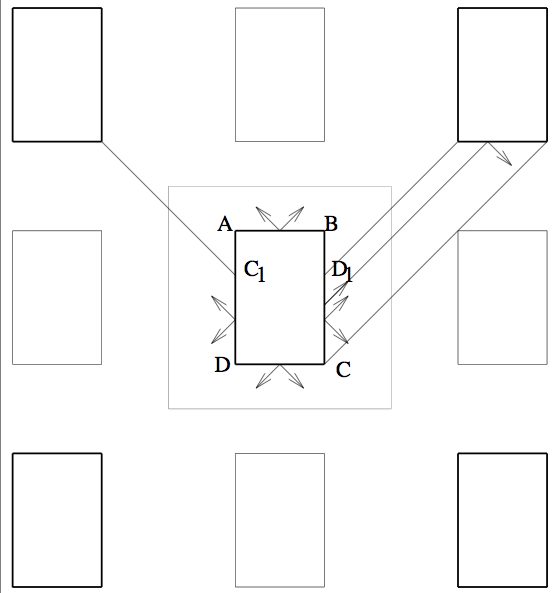}
\end{center}
\caption{\it The billiard table, angle $\pi /4$}
\label{bill-fig1}
\end{figure}

A natural question is the existence of a sequence of (deterministic) times along which the process $(X^t)$ after normalisation has a limit distribution.
For the second part of the next result, we need the mild assumption on the partial quotients of $\alpha$:
\begin{eqnarray}
n^{-\frac12} \,\sum_{j=1}^{[\ln n]} a_{j+1} \to 0 \label{mildhyp1}
\end{eqnarray}
This condition is satisfied by a.e. $\alpha$. Indeed, since $\alpha \to \left(a_1(\alpha)\right)^\frac12$ is integrable, 
we have for a.e. $\alpha$ by the ergodicity of the Gauss map,
$\sup_N  \, N^{-1} \, \left(\left(a_1(\alpha)\right)^\frac12 + ... + \left(a_N(\alpha)\right)^\frac12\right) < \infty$,
which implies $\sup_N \, N^{-2} \, (a_1(\alpha) + ... + a_N(\alpha)) < \infty$

\begin{thm} For almost all ratio $a/b$ of the lengths of the sides of the rectangular obstacles, there is a sequence $(L_n)$ such that
$S(L_n, \Psi)$ has asymptotically after normalization a non-degenerate 2-dimensional normal distribution
(with the uniform measure on a rectangular obstacle as initial distribution). Moreover under the assumption (\ref{mildhyp1}), also
$n^{-\frac12} X^{c L_n}$ converges in distribution to a non-degenerate normal law.
\end{thm}
\proof \  The first part of the theorem follows from Theorem \ref{vectCLT1}.

For the second part, observe, first, that $T_n(x) \sim c \, n$ by the ergodic theorem, where $c$ is defined in (\ref{centerpsi}).

Let us show that the distance between $X^{T_{L_n}(x)} (x) $ and $X^{c \, L_n} (x)$ is small compared with the norm of the vector $S(L_n, \Psi)(x)$,
which is of order $\sqrt n$. 

Denoting by $d$ the euclidean distance in the plane, the distance $d(X^{T_{L_n}(x)} (x) , X^{c \, L_n} (x))$ corresponds to a travel  
of the ball during a lapse of time $|T_{L_n}(x) - c \, L_n|$.

The (time) difference $T_{L_n}(x) - c \, L_n$ coincides with $\psi^0_{L_n}$, the ergodic sum of $\psi^0$ at time $L_n$.
The estimation of the variance along the sequence $L_n$ can be applied to the function $\psi^0$ which is centered and belongs to the class $\cal C$.

We get: $\|\psi^0_{L_n}\|_2^2 \leq C n$. Therefore, for $\varepsilon >0$, there is a constant $M= M(\varepsilon)$ such that
$|\psi^0_{L_n}(x)| \leq M \sqrt n$ on a set $B_\varepsilon$ of measure $\geq 1 - \varepsilon$.

If the ball starts from a point $x$ on the boundary of the obstacle located on some cell $p \in \Z^2$ and travels during a lapse of time $t$ such that 
$\psi_n(x) \leq t < \psi_{n+1}(x)$, then the ball hits $2n$ or $2n+1$ obstacles and reaches a cell at a (uniformly) bounded distance of the cell $p + S(n, \Psi)(x) \in \Z^2$.
By Ostrowski expansion, cf. (\ref{OstroBnd1}), we have $\|S(n, \Psi)(x)\| \leq C  \sum_1^{\ell+1} a_{j}$, if $q_\ell \leq n < q_{\ell+1}$.
By the ergodic theorem, $\psi_n(x) / n$ is of order $c$, so $n$ is of order $t/c$ and less than $t/ \inf_{x \in \T} \psi(x)$.

If we know only that $t \leq K$ for some $K$, then the distance between the starting point and the final point after time $t$ is 
$\leq C  \sum_1^{r+1} a_{j} + C'$, where $q_r \leq K/c< q_{r+1}$ and $C, C'$ are the constants $C= V(\psi), C'= \sup_{x \in \T} \psi(x)$.

This shows that, on the set $B_\varepsilon$, the distance between $X^{T_{L_n}(x)} (x)$ and $X^{c \, L_n} (x)$ is bounded by:
$C \sum_1^{s+1} a_{j} + C'$, where $s$ is such that $q_s \leq CM \sqrt n /c< q_{s+1}$.

In view of the (at least) exponential growth of the $q_n$'s, $s$ is at most of order $C \ln n$. We conclude by (\ref{mildhyp1}) that 
$n^{-\frac12} \, d(X^{T_{L_n}(x)} (x) , X^{c \, L_n} (x)) \to 0$ in probability. Therefore,
$n^{-\frac12} X^{T_{c L_n}}(x)$ has the same normal limit in distribution as $n^{-\frac12} \, X^{T_{L_n}(x)} (x) = n^{-\frac12} \, S(L_n, \Psi)(x)$.
\eop

\goodbreak
\section{\bf Appendix} \label{ASIP0}

\subsection{\bf CLT and ASIP for $\sum f_k(n_k.)$} \label{ASIPlacun}

\

In this appendix, we present an extension of a result of Berkes and Philipp \cite{BePh79} which was used in the proof of Theorem \ref{TCLsubseq0}.

Let $(n_k)$ be an increasing lacunary sequence of positive integers and $\rho := \inf_k {n_{k+1} \over n_k} > 1$.

For a fixed regular function $f$, the problem of the CLT for the sums $\sum_{j=1}^{N} f(n_j.)$ has been studied by several authors
(Zygmund, Salem, Kac, Fortet, then Gaposhkin, Berkes, Berkes and Philipp, ...). Two questions arise: non nullity of the variance, validity of
a CLT when the variance is non zero. The answer to the second question depends on arithmetic conditions on the sequence $(n_k)$, cf. Gaposhkin and
Subsection \ref{GapoApp}. For the first question, in general it is difficult to check the non nullity of the asymptotic variance
$\lim_N {1\over N} \|\sum_{k=1}^N f(n_k .)\|_2^2$, if it exists, except for special sequences or when Corollary \ref{rhoCond1} below can be used.

For $f \in L^2(\T)$ with Fourier coefficients $(c_j)$, we put $R(f, t) = (\sum_{|j| \geq t} |c_j|^2)^\frac12, \text{ for } t > 0$.
If $f$ is in $\cal C$, we have: $R(f, t) = (\sum_{|j| \geq t} |{\gamma_j \over j}|^2)^\frac12 \leq 2K(f) \, t^{-\frac12}$.

We introduce now several hypotheses.
\begin{hypo0} \label{hypofk} (H\ref{hypofk}) (regularity) We say that a sequence $\cal F = (f_k)_{k\geq 1}$ of real functions in $L^2(\T^1)$
with zero mean satisfies {$(H\ref{hypofk})$} if the following conditions hold:
\hfill \break 1) (uniform bound  for uniform and $L^2$-norms)
\begin{eqnarray}
&& M(\cal F) := \sup_k \|f_k\|_\infty < +\infty, \ \Phi= \Phi(\cal F) := \sup_k \|f_k\|_2 < +\infty;
\end{eqnarray}
2) (tail of the Fourier series) there is a finite constant $C_R= C_R(\cal F)$
and a constant $\gamma = \gamma(\cal F) > 0$ such that the tail of the Fourier series of $f_k$ satisfies uniformly in $k$:
\begin{eqnarray}
&&R(f_k, t) \leq R(t), \ \forall k, \text { with } R(t) \leq C_R \, t^{-\gamma}. \label{hypofkIneg}
\end{eqnarray}
\end{hypo0}
\begin{hypo0} \label{hyponk} (H\ref{hyponk}) (lacunarity of a sequence $(n_k)$) There is $\rho > 1$ such that
\begin{eqnarray}
n_{k+1} /n_k \geq \rho > 1, \ \forall k \geq 1. \label{(1.3)}
\end{eqnarray}
\end{hypo0}
\begin{hypo0}
\label{Gapo} (H\ref{Gapo}) (arithmetic condition)  For all integers $m \geq 1$, the following  holds for $\{n_k\}$:
\hfill \break Condition $(D_m)$: \ there is a constant $C$ such that the equation $t n_k \pm s n_\ell = \nu$, for $k>\ell$
and $t, s = 1, ..., m$, has at most $C$ solutions for any integer $\nu > 0$.
\end{hypo0}
The CLT for $\{f(n_k .)\}$ under (H\ref{hyponk}) and (H\ref{Gapo}) follows from Gaposhkin \cite{Ga70} for a given 
sufficiently smooth function $f$. Results of Berkes \cite{Be76} and Berkes and Philipp \cite{BePh79} give an approximation by a Wiener process:
\begin{thm} \cite{BePh79} \label{BerkPhilThm} Let $f$ be a 1-periodic Lipschitz centered function.
Assume that $(n_k)$ satisfies (H\ref{hyponk}) and  (H\ref{Gapo}). Assume moreover the condition
\begin{eqnarray}
&&\exists \, C >0,  N_0 \geq 1 \text{ such that }
\int_0^1 [\sum_{k=M+1}^{M+N} f(n_k x)]^2 \, dx \geq C N, \ \forall M \geq 0, \ \forall N \geq N_0. \label{varianceStrg}
\end{eqnarray}
Let $S_N = \sum_{k=1}^N f(n_k x)$. Then the sequence $(S_N, N \geq 1)$ can be redefined on a new
probability space (without changing its distribution) together with a Wiener process
$\zeta(t)$ such that
\begin{eqnarray*}
S_N = ~\zeta(\tau_N)+O(N^{1/2- \lambda}) \ \text{a.e.}
\end{eqnarray*}
where $\lambda >0$ is an absolute constant and $(\tau_N)$ is an increasing sequence of random
variables such that $\tau_N/\|S_N\|_2^2 \to 1$ a.s.
\end{thm}
We use the following slightly extended version in which the fixed Lipschitz function $f$ of Theorem \ref{BerkPhilThm} is replaced by a family $(f_k)$ satisfying a uniform 
boundedness condition and a uniform tail condition.  Moreover Condition (\ref{varianceStrg}) of Theorem \ref{BerkPhilThm} can be replaced by a weaker one.

\begin{thm} \label{ext-BerkPhil} Let $(f_k)$ satisfy (H\ref{hypofk}) and let $(n_k)$ be  a sequence of integers satisfying (H\ref{hyponk})
and (H\ref{Gapo}). Let $S_N(x) := \sum_{k=1}^N f_k(n_k x)$. Suppose that the following condition holds:
\begin{eqnarray}
&&\exists \, C>0,  N_0 \geq 1 \text{ such that }
\int_0^1 [S_N(x)]^2 \, dx \geq C N, \ \forall N \geq N_0. \label{minMN}
\end{eqnarray}
Then the process $(S_N)_{N \geq 1}$ can be redefined on a
probability space (without changing its distribution) together with a Wiener process $(\zeta(t))_{t \geq 0}$ such that
for an absolute constant $\lambda >0$ and an increasing sequence of random
variables $(\tau_N)$ satisfying $\tau_N/\|S_N\|_2^2 \to 1$ a.e. we have:
\begin{eqnarray}
S_N = \zeta(\tau_N) +o(N^{\frac12-\lambda}) \ \text{a.e.}
\end{eqnarray}
\end{thm}
By Corollary \ref{rhoCond1} below, for $\rho$ big enough, (\ref{minMN}) reduces to $\liminf_N \, \frac1N \sum_{k=1}^{N} \, \|f_k\|_2^2 > 0$. 

For the sake of conciseness, we do not reproduce the proof of this extended version which is an adaptation of the proofs in \cite{Be76} and \cite{BePh79}.

\vskip 3mm
{\it Quasi-orthogonality and variance}

\begin{lem} \label{lemA1} Let $f, g$ be in $L_0^2(\T^1)$, $\lambda_2 \geq \lambda_1$ two positive integers.  Then
\begin{eqnarray}
&&|\int_0^1 f(\lambda_1  x) \, \overline { g(\lambda_2 x) }\, dx| \leq R(f, {\lambda_2 \over \lambda_1}) \, \|g\|_2. \label{InegDec0}
\end{eqnarray}
\end{lem}
\proof \ (\ref{InegDec0}) follows from Parseval relation and
$$\int_0^1 f(\lambda_1  x) \, \overline { g(\lambda_2 x)} \, dx = \sum_{k, \ell: \, \lambda_1 k = \lambda_2 \ell} c_k(f) \,  \overline {c_{\ell}(g)}. \eop $$

\begin{lem} Let $\mathcal{F}=(f_k)$ satisfy (H\ref{hypofk}). If $(n_k)$ is a sequence of integers satisfying (\ref{(1.3)}) , then we have
$$\int_0^1\left(\sum_{k=1}^Nf_k(n_k x)\right)^2dx=\sum_{k=1}^N\| f_k\|_2^2+NB_N,$$
with $B_N\leq C\Phi(\mathcal{F})\rho^{-\gamma}$, where the constant $\gamma$ is the one given in (H\ref{hypofk}), $\rho$ is the constant of lacunarity in
(\ref{(1.3)}) and $C$ depends on the constant $C_R(\mathcal{F})$ defined in (\ref{hypofkIneg}). In particular
$$\int_0^1\left(\sum_{k=1}^Nf_k(n_k x)\right)^2dx=O(N).$$
\end{lem}
\Proof \  Putting $W_k = \sum_{\ell = 1}^{N-k} \int_0^{1} f_\ell(n_\ell x) \, f_{\ell+k}(n_{\ell + k} x) \, dx$, we have
\begin{eqnarray*}
&&\int_0^{1} (f_1(n_1 x) + ... + f_N(n_N x))^2  dx = \sum_{k = 1}^N \int_0^1 f_k^2(n_k x) dx +2 W_1 + ... + 2W_{N-1}.
\end{eqnarray*}
We have $\displaystyle {n_{\ell + 1} \over  n_\ell} \geq \rho > 1$ and by Lemma \ref{lemA1}, for $1 \leq k \leq N-1$:
$|W_k| \leq C_R(\cal F) \, \Phi(\cal F) \, N \, \rho^{-k\gamma}$.

It follows: $\displaystyle |W_1 + ... + W_{N-1}| \leq C_R(\cal F) \, \Phi(\cal F) \, N \, {\rho^{-\gamma} \over 1 - \rho^{-\gamma}}$.
\eop

\begin{cor} \label{rhoCond1} There is $\rho_0 > 1$ and $c >0$ depending on $\Phi(\Cal F), C_R(\Cal F), \gamma(\Cal F)$ such that, for $\rho \geq \rho_0$,
\begin{eqnarray}
\int_0^1 [\sum_{k=1}^{N} f_k(n_k x)]^2 \, dx \geq c \sum_{k=1}^{N} \, \|f_k\|_2^2, \  \forall N \geq N_0.
\end{eqnarray}
\end{cor}

A similar proof shows the following lemma.
\begin{lem} \label{varianceoN}
If $\mathcal{F}=(f_k)$ satisfies (H\ref{hypofk}) and  $(m_k)$ is a superlacunary sequence of positive integers (i.e., $n_{k+1}/n_k\rightarrow \infty$), then we have
$\int_0^1 \, (\sum_{k=1}^Nf_k(n_k x))^2 dx = \sum_{k=1}^N\| f_k\|_2^2+o(N)$.
\end{lem}

\vskip 3mm
\subsection{\bf A remark on a result of Gaposhkin} \label{GapoApp}

\

Let $(n_k)$ be a lacunary sequence of integers and $f$ a 1-periodic real function with some regularity. If the quotients $n_{k+1}/ n_k$ are integers,
the  central limit theorem (CLT) holds for the sums $\sum_{k=0}^{n-1}  f (n_k x)$. But Erd\"os and Fortet gave the example of $n_k = 2^k -1$
for which the function $f_0(x) := \cos(2\pi x) +\cos(4\pi x)$ does not satisfy the CLT.

Let us recall this counter-example (cf.~\cite{AiBe08}, \cite{CoLe11}).
If $f_0$ is as above and $Z_n(x) = {1\over \sqrt{n}} \sum_{k=1}^{n}  f_0(2^k x -x)$, then
\begin{eqnarray}
\mu\{x :  {1\over \sqrt{n}} Z_n(x) \leq t \} \rightarrow {1\over \sqrt{2\pi}} \int_0^1 (\int_{-\infty}^{t/ \cos y} e^{-t^2/2} dt)\  dy.\label{melange}
\end{eqnarray}
For the proof, observe that the sum $\sum_1^{n}[\cos(2\pi (2^k -1)x) + \cos(4\pi (2^k -1)x)]$ reads
\begin{eqnarray*}
&\cos(2\pi x) + \cos(2\pi (2^{n+1} -2)x) + \sum_2^{n} [\cos(2\pi (2^{k} -1)x) + \cos(2\pi (2^{k} -2)x)]\\
&=\cos(2\pi x) + \cos(2\pi (2^{n+1} -2)x) + 2 \cos(\pi x) \sum_2^{n} \cos(2\pi (2^k  - 3/2) x).
\end{eqnarray*}
The convergence then follows from the CLT in Salem and Zygmund \cite{SaZy48} and from the following lemma (see for example \cite{CoLe11}).
\begin{lem} \label{phi} Let $(Y_n)$ be a sequence of random variables defined on $([0,1], \PP)$
and $\cal{L}$ a distribution on $\R$ with characteristic function $\Phi$. The following conditions are  equivalent:
\hfill \break a) for every probability density $\psi$, the sequence $(Y_n)$ converges in distribution to
$\cal{L}$ under the measure $\psi \PP$; \hfill \break b) for every interval $I$,
${1\over \mu(I)} \mu\{x \in I: Y_n(x) \leq t) \rightarrow {\cal L}(]-\infty, t])$,
$\forall t \in \R$; \hfill \break c) for every Riemann integrable function $\varphi$, the sequence $(\varphi Y_n)$ converges in distribution
to a limit with characteristic function $\int_0^1 \Phi(\varphi(y)\,t) \ dy$. In particular, if ${\cal L} = {\cal N}(0,1)$, the sequence
$(\varphi Y_n)$ converges in distribution to a limit with characteristic function $\int_0^1 e^{-{1\over 2}\varphi^2(y)\, t^2} \, dy$,
a mixture of Gaussian distributions.
\end{lem}

\goodbreak
{\it Description of a result of Gaposhkin}

Gaposhkin has introduced an arithmetical condition on $(n_k)$ so that the CLT should be true. Actually, he has given an answer to a slightly
different problem. For simplicity, we consider only trigonometric polynomials.

Let $(\lambda_{M,k}, \ 1 \leq k \leq M, M \geq 1)$ be an array of non negative numbers. We say that Property $(P)$ holds if
for every  $M$, $\sum_{k=1}^M\lambda_{M,k}^2 =  1$ and $\lim_M\max_k\|\lambda_{M,k}\|=0$.

\begin{thm} (Gaposhkin) \label{gapo} Let $(n_k)$ be a lacunary sequence satisfying the arithmetic condition (H\ref{Gapo})
(i.e., $(D_m)$ for every $m$). Then the following strong version of the CLT holds for every trigonometric polynomial $f$:
if $(\lambda_{M,k})$ is an array with property $(P)$, for every measurable subset $E$ of $[0,1]$ with positive measure,
we have, with $\lambda_{M}( f )^2 :=\E\left((\sum_{k=1}^\infty\lambda_{M,k} f (n_k\cdot))^2\right)$,
\begin{eqnarray}
&&\frac{1}{\PP(E)}\PP (x\in E\ :\ \lambda_{M}( f )^{-1}\sum_k \lambda_{M,k} f (n_k x)<y)
\longrightarrow_{M\rightarrow \infty} \frac{1}{\sqrt{2\pi}}\int_{-\infty}^ye^{-t^2/2}dy.
\end{eqnarray}
\end{thm}
A possible choice for $(\lambda_{M,k})$ is as in the "classical" version of the CLT:
\begin{eqnarray}
\lambda_{M,k}=\frac{1}{\sqrt{M}}\ {\rm if}\ 1\leq k\leq M,\ \ = 0\ {\rm otherwise.}
\end{eqnarray}
Gaposhkin also has shown that if $(D_m)$ is not satisfied for every $m \geq 1$, one can find a trigonometric polynomial $ f $ for which the above strong
version of the result is not true anymore. But this does not mean that classical CLT is not true, as we will see.

For an integer $a>1$, let us consider the following subset of $\N^*$:
$I_a=\cup_{n\geq 1}\{k\in \N\ /\ \ k \in [n^a,n^a+n]\}$. Let $(n_k)$ be the sequence
$n_k=2^k\ {\rm if}\ k\notin I_a ,\ n_k=2^k-1\ {\rm if} \ k \in I_a$.

This sequence does not satisfy the condition (H\ref{Gapo}). It is easy to find a family $\lambda_{M,k}$ for which the conclusion of the preceding theorem
is not true. It suffices to consider the family $\lambda_{M,k}=0\ {\rm if}\ k \notin I_a$,
$\lambda_{M,k}=\frac{1}{\sqrt{M}}\ {\rm if}\ k\ {\rm is\ one\ of\ the\ first} \ M \ {\rm elements \ of}\ I_a$.

Following Kac, Fortet, for this choice of $\lambda_{M,k}$ and for $f_0$,
the central limit theorem is not satisfied (this is what Gaposhkin did to show that his condition is necessary).
For $M=\frac{N(N+1)}{2}$, $\lambda_{M,k}$ is either 0 or $(\frac{N(N+1)}{2})^{-\frac12}$. The quantity
$(\frac{N(N+1)}{2})^\frac12 \,\sum_{k=1}^{+\infty}\lambda_{\frac{N(N+1)}{2},k} \, f (n_k x)$
$=(\frac{N(N+1)}{2})^\frac12 \,\sum_{k=1}^{N^a+N}\lambda_{\frac{N(N+1)}{2},k} \, f (n_k x)$ reads after computation:
$\sum_{j=1}^N\cos(2\pi(2^{j^a}-1) x)+\sum_{j=1}^N \cos(2\pi (2^{j^a+j+1}-2)x)$ $+\cos(\pi x)\sum_{j=1}^N (\sum_{k=j^a+1}^{j^a+j}
\cos(2\pi(2^k - \frac32) x))$. The first sums have $N$ terms, the third one $\frac{N(N+1)}{2}$.

By the CLT in \cite{SaZy48}, $(\frac{N(N+1)}{2})^{-\frac12} \, \sum_{j=1}^N \, [\sum_{k=j^a+1}^{j^a+j}\cos(2\pi(2^k-3/2) x)]$
converges toward a gaussian variable. Lemma \ref{phi} then implies that
$$(\frac{N(N+1)}{2})^{-\frac12}\cos(\pi x)\sum_{j=1}^N \, [\sum_{k=j^a+1}^{j^a+j}\cos(2\pi(2^k-3/2) x)]$$
converges toward a distribution similar to the one appearing in (\ref{melange}). The same is true for
$\sum_{k=1}^{N^a+N}\lambda_{\frac{N(N+1)}{2},k} f (n_k x)$ as
$(\frac{N(N+1)}{2})^{-\frac12} \, \sum_{j=1}^N\left(\cos(2\pi(2^{j^a}-1) x)+\cos(2\pi (2^{j^a+j+1}-2)x)\right)$
converges to 0. One easily shows that $\sum_{k=1}^{\infty}\lambda_{M,k} f (n_k x)$ has the same limit when $M \to \infty$.

We now consider another choice for $\lambda_{M,k}$, the classical choice $\lambda_{M,k}=1/\sqrt{M}$ if $k=1,\ldots,M$, 0 otherwise.
Let us begin by a simple remark. Let  $f$ be a trigonometric polynomial and $(n_k)$ and $(q'_k)$ two sequences such that
$\#\{k\in\{1,\ldots,n\}: n_k\neq q'_k\}=o(\sqrt{n})$. Suppose that $n^{-1/2}\sum_{k=1}^n f (n_k\cdot)$ converges in law toward a random variable $Y$.

Then, $n^{-1/2}\sum_{k=1}^n f (q'_k\cdot)$ also converges toward $Y$, as $\sum_{k=1}^n f (n_k\cdot)-\sum_{k=1}^n f (q'_k\cdot)=o(\sqrt{n})$.
Let us take, for $(n_k)$ the sequence $(2^k)_{k\geq 1}$, and for $(q'_k)$ : $q'_k=2^k\ {\rm if}\ k\notin I_a ,\ q'_k=2^k-1\ {\rm if} \ k\in I_a$. Then
$$\#\{k\in\{1,\ldots,n\}: \ n_k\neq q'_k\}=\#\left(\{1,\ldots,n\}\cap I_a\right)=O\left(n^{2/a}\right).$$

From this, we deduce that for $a\geq 5$ and $ f $ a trigonometric polynomial not of the form $\psi(\cdot)-\psi(2\cdot)$, both sums
$n^{-1/2}\sum_{k=1}^n f (n_k\cdot)\ {\rm and}\ n^{-1/2}\sum_{k=1}^n f (q'_k\cdot)$
converge to the same non-degenerate Gaussian law (recall that CLT is satisfied for the sequence $(2^k)_{k\geq 1}$).

In other words the necessary and sufficient condition of Gaposhkin is only necessary for the strong result, i.e., Theorem \ref{gapo}.
The classical version can be true for examples without this condition holding. This is not a surprise:
the above remark shows that a sufficiently rare modification of the sequence $(n_k)$ can not be seen anymore at infinity after normalisation
by $\sqrt{n}$. The arithmetic condition (H\ref{Gapo}) is much more rigid.

\bibliographystyle{amsalpha}

\providecommand{\bysame}{\leavevmode\hbox to3em{\hrulefill}\thinspace}
\providecommand{\MR}{\relax\ifhmode\unskip\space\fi MR }
\providecommand{\MRhref}[2]{%
  \href{http://www.ams.org/mathscinet-getitem?mr=#1}{#2}
}
\providecommand{\href}[2]{#2}
\begin{thebibliography}{}

\end{thebibliography}


\begin{thebibliography}{12345678}

\bibitem[AiBe08]{AiBe08}
Aistleitner (C.), Berkes (I.): On the central limit theorem for
$f(n_k x)$, {Probab. Theory Related Fields}, {146}, {2010}, {1-2}, p. {267-289}.

\vskip 1mm
\bibitem [Be76]{Be76}
Berkes (I.):  On the asymptotic behaviour of $\sum f(n_kx)$: I. Main theorems, II. Applications, Z. Wahrscheinlichkeitstheorie
verw. Gebiete, 34 (1976), 319-345 and 347-365.

\vskip 1mm
\bibitem [BePh79]{BePh79}
Berkes (I.), Philipp (W.): An a.e. invariance principle for lacunary series $f(n_k x)$,
Acta Math. Acad. Sci. Hungar. 34 (1979), no. 1-2, 141-155.

\vskip 1mm
\bibitem[CoGu12]{CoGu12}
Conze (J.-P.), Gutkin (E.): On recurrence and ergodicity for geodesic flows on non-compact periodic polygonal surfaces.
Ergodic Theory Dynam. Systems, 32(2): 491-515, 2012.

\vskip 1mm
\bibitem[CoLe11]{CoLe11}
Conze (J.-P.), Le Borgne (S.): Limit law for some modified ergodic sums. Stoch. Dyn., 11(1): 107-133, 2011.

\vskip 1mm
\bibitem [Ga70]{Ga70}
Gaposhkin (V.F.):  On the central limit theorem for some weakly dependent sequences (in
Russian), Teor. Verojatn. i Primenen, 15 (1970).

\vskip 1mm
\bibitem[Go10]{Go10}
Gou\"ezel (S.), Almost sure invariance principle for dynamical systems
by spectral methods, Ann. Probab. 38 (2010), no. 4, 1639-1671.

\vskip 1mm
\bibitem[HW80] {HW80}
Hardy (J.) and Weber (J.): Diffusion in a periodic wind-tree model, J. Math. Phys. {\bf 21} (1980), 1802-1808.

\vskip 1mm
\bibitem[Hu09]{Hu09}
Huveneers (F.): Subdiffusive behavior generated by irrational rotations, Ergodic Theory Dynam. Systems 29 (2009), no. 4, 1217-1233.

\vskip 1mm
\bibitem[Is06]{Is06}
Isola (S.): Dispersion properties of ergodic translations, Int. J. Math. Math. Sci. 2006, Art. ID 20568, 20 pp.

\vskip 1mm
\bibitem [KSZ48]{KSZ48}
Kac (M.), Salem (R.), Zygmund (A.): A gap theorem. Trans. Amer. Math. Soc. 63, 235-243 (1948).

\vskip 1mm
\bibitem [Kh37] {Kh37}
Khinchin (A.Ya.): {\it Continued fractions}, Dover Publications, Mineola, N.Y, 1997.

\vskip 1mm
\bibitem [PhSt75]{PhSt75}
W. Philipp (W.), Stout (W. F.): Almost sure invariance principles for partial sums of weakly dependent random variables.
Mem. Amer. Math. Soc. 2 (1975), 2, no. 161.

\vskip 1mm
\bibitem[SaZy48]{SaZy48}
Salem (R.), Zygmund (A.):  On lacunary trigonometric series. II. Proc. Nat. Acad. Sci. U. S. A. 34, (1948) p. 54-62.
\end{thebibliography}

\end{document}